\documentclass[reqno, 12pt]{amsart}
\usepackage{a4wide}
\usepackage{amscd}
\input xy
\xyoption{all}

\newtheorem{thm}{Theorem}[section]

\newtheorem{prop}[thm]{Proposition}

\newtheorem{lm}[thm]{Lemma}

\numberwithin{equation}{section}

\font\fr=eufm10 scaled \magstep1

\def\R{\mathbb R}

\def\cv{\hbox{\fr v}}

\def\ad{\hbox{\rm ad}}

\def\id{\hbox{\rm id}}

\begin{document}

\title{Cohomogeneity one actions on symmetric spaces of noncompact type}
\author{\textsc{J\"{u}rgen Berndt and Hiroshi Tamaru}}
\address{Department of Mathematics \\ King's College London \\ Strand \\ London \\ WC2R 2LS \\ United Kingdom}
\email{jurgen.berndt@kcl.ac.uk}
\address{Department of Mathematics \\ Graduate School of Science\\ Hiroshima University \\
1-3-1 Kagamiyama \\ Higashi-Hiroshima, 739-8526\\ Japan}
\email{tamaru@math.sci.hiroshima-u.ac.jp}

\begin{abstract}
An isometric action of a Lie group on a Riemannian manifold is of cohomogeneity one if
the corresponding orbit space is one-dimensional. In this article we develop
a conceptual approach to the classification
of cohomogeneity one actions on Riemannian symmetric spaces of noncompact type
in terms of orbit equivalence.
As a consequence, we find many new examples of cohomogeneity one actions
on Riemannian symmetric spaces of noncompact type. We apply our conceptual approach
to derive explicit classifications of cohomogeneity one actions on some symmetric spaces.
\end{abstract}

\maketitle
\thispagestyle{empty}



\section{Introduction}

The cohomogeneity of an isometric action on a Riemannian manifold is the rank of the normal
bundle of a principal orbit of the action.
Thus, for a cohomogeneity one action, the principal orbits are hypersurfaces.
Cohomogeneity one actions have been of much recent interest
in the context of constructing geometric structures on manifolds.

The main focus in this article is on the classification of such actions on Riemannian symmetric spaces. A remarkable result by Hsiang and Lawson (\cite{HL}) states that every cohomogeneity one action on the round sphere $S^n$ is orbit equivalent to the action on the sphere $S^n$ which is induced from the isotropy representation of an $(n+1)$-dimensional Riemannian symmetric space of rank two. Cohomogeneity one actions on the other compact simply connected Riemannian symmetric spaces of rank one, that is, the projective spaces over the normed real division algebras ${\mathbb C}$, ${\mathbb H}$ and ${\mathbb O}$, were obtained by Takagi (\cite{Tak}) and Iwata (\cite{Iw1},\cite{Iw2}). Kollross (\cite{Ko}) derived the classification of cohomogeneity one actions on irreducible compact simply connected Riemannian symmetric spaces of higher rank. The classification for the reducible case is still outstanding.

In the noncompact case one needs to develop different techniques due to the noncompactness of the isometry groups. This can already be seen when considering cohomogeneity one actions on the Euclidean space ${\mathbb E}^n$. A group theoretical approach as in the compact case leads immediately to difficulties. However, there is a simple geometric solution to the problem. A principal orbit of a cohomogeneity one action is a hypersurface with constant principal curvatures, also known as an isoparametric hypersurface. Isoparametric hypersurfaces in Euclidean spaces were classified by Somigliana (\cite{Som}), Levi-Civita (\cite{Lev}) and Segre (\cite{Seg}), and it is easy to verify from their results that all complete isoparametric hypersurfaces in Euclidean spaces are homogeneous and hence principal orbits of cohomogeneity one actions.
A similar approach leads to the classification of cohomogeneity one actions on real hyperbolic spaces by using the classification of isoparametric hypersurfaces in real hyperbolic spaces  by Cartan (\cite{Car}). However, this approach is successful only in these two cases. For example, the classification of hypersurfaces with constant principal curvatures in complex hyperbolic spaces is not yet known. In this article we present a conceptual approach for classifying cohomogeneity one actions on Riemannian symmetric spaces of noncompact type up to orbit equivalence. Two actions are said to be orbit equivalent if there exists an isometry of the space mapping the orbits of one action onto the orbits of the other action.

The orbit space of a cohomogeneity one action of a connected Lie group on a connected complete Riemannian manifold $M$ is homeomorphic to the closed bounded interval $[0,1]$, the closed unbounded interval $[0,\infty)$, the circle $S^1$ or the real line ${\mathbb R}$, each of them equipped with their standard topology. If $M$ is a Riemannian symmetric space of noncompact type, then for topological reasons the orbit space must be homeomorphic to either ${\mathbb R}$ or $[0,\infty)$. In the first case the orbits form a Riemannian foliation on $M$, and in the second case there is exactly one singular orbit and the principal orbits are the tubes around this singular orbit.

Let $M = G/K$ be a connected Riemannian symmetric space of noncompact type
and $r = {\rm rank}(M)$, where $G$ is the identity component of the isometry
group of $M$ and $K$ is the isotropy subgroup of $G$ at a point $o \in M$.
Let $H$ be a connected subgroup of $G$ which acts on $M$ with cohomogeneity one. The case when the orbits of $H$ form a Riemannian foliation on $M$ has been dealt with by the authors in \cite{BT1} for irreducible symmetric spaces $M$. We therefore assume that the action has a singular orbit $W$. Without loss of generality we may assume that $o \in W$.
The subgroup $H$ is contained in a connected maximal proper subgroup $L$ of $G$. It follows from work by Mostow (\cite{Mo2}) that $L$ is either reductive or the identity component of a parabolic
subgroup of $G$. For the reductive case we show in Theorem \ref{maximal} that $H$ and $L$ are orbit equivalent and that the singular orbit $W$ is a totally geodesic submanifold in $M$. We now assume that $L$ is the connected identity component of a parabolic subgroup of $G$.

The conjugacy classes of parabolic subgroups of $G$ can be parametrized by
the subsets $\Phi$ of a set $\Lambda = \{\alpha_1,\ldots,\alpha_r\}$
of simple roots of a restricted root system of the semisimple Lie algebra
${\mathfrak g}$ of $G$.
The maximal proper parabolic subgroups correspond to subsets $\Phi$ of $\Lambda$ with cardinality $|\Phi|$ equal to $r-1$. For $\Phi = \emptyset$ we obtain a minimal parabolic subgroup of $G$.
Let $Q_\Phi$ be the parabolic subgroup of $G$ associated with the subset
$\Phi$ of $\Lambda$. We construct new examples of cohomogeneity one actions on $M$ from the
Langlands decomposition and from the Chevalley decomposition of $Q_\Phi$.

The Langlands decomposition is of the form
$Q_\Phi = M_\Phi A_\Phi N_\Phi$, where $M_\Phi$ is reductive, $A_\Phi$ is abelian and
$N_\Phi$ is nilpotent. The orbit $B_\Phi = M_\Phi \cdot o$ is a semisimple Riemannian symmetric space of noncompact type with ${\rm rank}(B_\Phi) =|\Phi|$, unless $\Phi = \emptyset$
in which case the orbit consists just of the point $o$. The symmetric space $B_\Phi$ is
embedded totally geodesically in $M$ and is also known as a boundary component of $M$ as it
arises naturally in the maximal Satake compactification of $M$.
The orbit $A_\Phi \cdot o$
is a Euclidean space ${\mathbb E}^{r - |\Phi|}$ of dimension $r - |\Phi|$
embedded in $M$ as a totally geodesic submanifold.
If $H_\Phi$ is a connected
subgroup of the isometry group of $B_\Phi$ acting on $B_\Phi$ with cohomogeneity one,
then $H = H_\Phi A_\Phi N_\Phi$ is a connected subgroup of $Q_\Phi \subset G$ acting on $M$
with cohomogeneity one. We call this the canonical extension of the cohomogeneity one action
on the boundary component $B_\Phi$ to the symmetric space $M$.

The Chevalley decomposition is of the form $Q_\Phi = L_\Phi N_\Phi$, where $L_\Phi = M_\Phi A_\Phi$
is reductive.
The orbit $F_\Phi = L_\Phi \cdot o$ is isometric to the Riemannian product $B_\Phi \times {\mathbb E}^{r - |\Phi|}$
and embedded in $M$ as a totally geodesic submanifold.
Let ${\mathfrak n}_\Phi$ be the Lie algebra of $N_\Phi$, and denote by $H^\Phi$ the sum
of the dual root vectors of the simple roots in $\Lambda \setminus \Phi$. The vector $H^\Phi$ induces
a gradation $\bigoplus_{\nu \geq 1} {\mathfrak n}_\Phi^\nu$ of ${\mathfrak n}_\Phi$
by defining ${\mathfrak n}_\Phi^\nu$ as the sum of all root
spaces corresponding to positive roots $\alpha$ with $\alpha(H^\Phi) = \nu \geq 1$.
Let ${\mathfrak v}$ be a subspace of ${\mathfrak n}_\Phi^1$ with dimension $\geq 2$.
Then ${\mathfrak n}_{\Phi,{\mathfrak v}} = {\mathfrak n}_\Phi \ominus {\mathfrak v}$ is a
subalgebra of ${\mathfrak n}_\Phi$. Denote by $N_{\Phi,{\mathfrak v}}$ the corresponding
connected subgroup of $N_\Phi$. Assume that
the normalizer $N_{L_\Phi}({\mathfrak n}_{\Phi,{\mathfrak v}})$ of
${\mathfrak n}_{\Phi,{\mathfrak v}}$ in $L_\Phi$ acts transitively on $F_\Phi$ and
that the normalizer $N_{K_\Phi}({\mathfrak n}_{\Phi,{\mathfrak v}})$ of
${\mathfrak n}_{\Phi,{\mathfrak v}}$ in $K_\Phi = L_\Phi \cap K$ acts
transitively on the unit sphere in ${\mathfrak v}$. Note that
$N_{K_\Phi}({\mathfrak n}_{\Phi,{\mathfrak v}})$ coincides with the normalizer
$N_{K_\Phi}({\mathfrak v})$ of ${\mathfrak v}$ in $K_\Phi$.
 Then
$H_{\Phi,{\mathfrak v}} = N^o_{L_\Phi}({\mathfrak n}_{\Phi,{\mathfrak v}})N_{\Phi,{\mathfrak v}}$
is a connected subgroup of $Q_\Phi = L_\Phi N_\Phi$ which acts on $M$ with cohomogeneity one
and singular orbit $H_{\Phi,{\mathfrak v}} \cdot o$. We provide some explicit examples of such actions below.

We put ${\mathfrak a} = {\mathfrak a}_\emptyset$ and ${\mathfrak n} = {\mathfrak n}_\emptyset$. Then ${\mathfrak g} = {\mathfrak k} \oplus {\mathfrak a} \oplus {\mathfrak n}$ is an Iwasawa decomposition of ${\mathfrak g}$, and the connected solvable subgroup $AN$ of $G$ with Lie algebra ${\mathfrak a} \oplus {\mathfrak n}$ acts simply transitively on $M$. Therefore $M$ is isometric to $AN$ equipped with a suitable left-invariant Riemannian metric.

Let $\ell$ be a one-dimensional linear subspace of ${\mathfrak a}$. Then ${\mathfrak h}_\ell =  ({\mathfrak a} \ominus \ell) \oplus {\mathfrak n}$ is a codimension one subalgebra of ${\mathfrak a} \oplus {\mathfrak n}$, and hence the connected subgroup $H_\ell$ of $G$ with Lie algebra ${\mathfrak h}_\ell$ acts on $M$ with cohomogeneity one. The orbits form a Riemannian foliation on $M$ whose orbits are pairwise isometrically congruent.

Let $\ell$ be a one-dimensional linear subspace of a simple root space ${\mathfrak g}_{\alpha_i}$. Then ${\mathfrak a} \oplus ({\mathfrak n} \ominus \ell)$ is a codimension one subalgebra of ${\mathfrak a} \oplus {\mathfrak n}$, and hence the connected subgroup of $G$ with Lie algebra ${\mathfrak a} \oplus ({\mathfrak n} \ominus \ell)$ acts on $M$ with cohomogeneity one. The orbits of this action form a Riemannian foliation on $M$, and there is exactly one minimal orbit. Moreover, assume that $\ell$ and $\ell^\prime$ are two one-dimensional linear subspaces of ${\mathfrak g}_{\alpha_i}$. Then the cohomogeneity one actions induced from ${\mathfrak a} \oplus ({\mathfrak n} \ominus \ell)$ and ${\mathfrak a} \oplus ({\mathfrak n} \ominus \ell^\prime)$ are orbit equivalent. Therefore, for each choice of simple root $\alpha_i \in \Lambda$ we get exactly one cohomogeneity one action up to orbit equivalence. We denote by $H_i$ one of the connected subgroups of $G$ constructed in this manner.

We can now formulate the main result of this article.

\begin{thm} \label{maintheorem}
Let $M = G/K$ be a connected irreducible Riemannian symmetric space of noncompact type and with rank $r$,
and let $H$ be a connected subgroup of $G$ which acts on $M$ with cohomogeneity one. Then either statement {\rm (1)} or statement {\rm (2)} holds:
\begin{itemize}
\item[(1)] The orbits form a Riemannian foliation on M and one of the following two cases holds:
\begin{itemize}
\item[(i)] All orbits are isometrically congruent to each other, and the action of
$H$ is orbit equivalent to the action of $H_\ell$ for some one-dimensional linear
subspace $\ell$ of ${\mathfrak a}$.
\item[(ii)] There exists exactly one minimal orbit, and the action of $H$ is orbit equivalent the action of $H_i$ for some $i \in \{1,\ldots,r\}$.
\end{itemize}
\item[(2)] There exists exactly one singular orbit and one of the following two cases holds:
\begin{itemize}
\item[(i)] $H$ is contained in a maximal proper reductive subgroup $L$ of $G$, the actions of $H$ and $L$ are orbit equivalent, and the singular orbit is totally geodesic in $M$.
\item[(ii)] $H$ is contained in a maximal proper parabolic subgroup $Q_\Phi$ of $G$ and one of the following two cases holds:
\begin{itemize}
\item[(a)] The action of $H$ is orbit equivalent to the canonical extension of a cohomogeneity one action
with a singular orbit on the boundary component $B_\Phi$ of $M$.
\item[(b)] The action of $H$ is orbit equivalent to a cohomogeneity one action on $M$ given by $H_{\Phi,{\mathfrak v}}$ for some subspace ${\mathfrak v} \subset {\mathfrak n}_\Phi^1$ with $\dim {\mathfrak v} \geq 2$.
\end{itemize}
\end{itemize}
\end{itemize}
\end{thm}

{\it Remarks.} 1. Consider the Dynkin diagram associated to the simple roots $\Lambda$. Each  symmetry $\sigma$ of the Dynkin diagram gives rise to an automorphism $F_\sigma$ of ${\mathfrak a}$.

In case (1)(i), assume that $\ell$ and $\ell^\prime$ are two one-dimensional linear subspaces of ${\mathfrak a}$. Then the cohomogeneity one actions induced from $({\mathfrak a} \ominus \ell) \oplus {\mathfrak n}$ and $({\mathfrak a} \ominus \ell^\prime) \oplus {\mathfrak n}$ are orbit equivalent if and only if there exists a Dynkin diagram symmetry $\sigma$ such that $F_\sigma(\ell) = \ell^\prime$. The cohomogeneity one actions of type (1)(i) are therefore parametrized by ${\mathbb R}P^{r-1}/{\mathfrak S}$, where ${\mathbb R}P^{r-1}$ is the real projective space of the real vector space ${\mathfrak a}$ and ${\mathfrak S}$ is the finite group of automorphisms of ${\mathbb R}P^{r-1}$ which is induced by the automorphisms $F_\sigma$ of ${\mathfrak a}$. For details we refer to Theorem 3.5 in \cite{BT1}.

In case (1)(ii), let $i,j \in \{1,\ldots,r\}$. The actions of $H_i$ and $H_j$ are orbit equivalent if and only if there exists a Dynkin diagram symmetry $\sigma$ such that $\sigma(\alpha_i) = \alpha_j$. The cohomogeneity one actions of type (1)(ii) are therefore parametrized by $\{1,\ldots,r\}/{\mathfrak S}$, where ${\mathfrak S}$ is the finite group of automorphisms of the Dynkin diagram. For details we refer to Theorem 4.8 in \cite{BT1}.

2. There is a well-known concept of duality between Riemannian symmetric spaces of noncompact type and Riemannian symmetric space of compact type. A totally geodesic submanifold $W$ of $M$ corresponds via this duality to a totally geodesic submanifold $W^*$ in the dual Riemannian symmetric space $M^*$ of compact type. A cohomogeneity one action of $H$ on $M$ with a totally geodesic singular orbit $W$ then gives rise to a cohomogeneity one action on $M^*$ of some connected subgroup $H^*$ of the isometry group of $M^*$. Using the classification by Kollross (\cite{Ko}) of cohomogeneity one actions on irreducible Riemannian symmetric spaces of compact type, and the concept of reflective submanifolds, the authors determined in \cite{BT2} all totally geodesic submanifolds in irreducible Riemannian symmetric spaces of noncompact type which arise as a singular orbit of a cohomogeneity one action. There are exactly five totally geodesic submanifolds which are not reflective, and mysteriously these are all related to the exceptional Lie group $G_2$. We refer to \cite{BT2} for further details.

We point out here that the explicit classification of totally geodesic submanifolds in reducible Riemannian symmetric spaces of noncompact type which arise as a singular orbit of a cohomogeneity one action is still an open problem.

3. The concept of canonical extension in (2)(ii)(a) suggests of course a rank reduction for the classification. However, since the boundary component $B_\Phi$ can be reducible, we encounter the same problem we discussed at the end of the previous remark.

4. We do not have an explicit classification of the groups $H_{\Phi,{\mathfrak v}}$ arising in (2)(ii)(b). However, our calculations indicate that there are only few examples which cannot be constructed via (2)(i) or (2)(ii)(a). The first author and Br\"{u}ck constructed in \cite{BB} new examples on the hyperbolic spaces over the normed real division algebras ${\mathbb C}$, ${\mathbb H}$ and ${\mathbb O}$. The authors proved in \cite{BT3} that there are no further examples in the cases of ${\mathbb C}$ and ${\mathbb O}$, but for ${\mathbb H}$ the problem remains open. In this article we construct two new cohomogeneity one actions with this method, one on $G_2^2/SO_4$ and one on $G_2^{\mathbb C}/G_2$. Although we checked many other symmetric spaces, we could not find any further examples and start to believe that there are none apart from the obvious ones on reducible symmetric spaces obtained from the known examples on irreducible symmetric spaces.

This article is organised as follows. In Section 2 we outline basic material about parabolic subalgebras of semisimple real Lie algebras, and relate this to the geometry of Riemannian symmetric spaces of noncompact type. In Section 3 we show first that a proper maximal reductive subgroup of the isometry group of a Riemannian symmetric space of noncompact type cannot act transitively on the space. We then relate cohomogeneity one actions to actions of reductive and parabolic subgroups. In Section 4 we present two new methods for constructing cohomogeneity one actions with a singular orbit on Riemannian symmetric spaces of noncompact type. In Section 5 we proof the main result of this article. In Section 6 we apply the main result to derive explicit classifications of cohomogeneity one actions on some Riemannian symmetric spaces of noncompact type and rank $2$.

\section{Parabolic subalgebras}
\label{parabolic}

In this section we recall the construction of the parabolic
subalgebras of real semisimple Lie algebras (see e.g.\ \cite{BJ},
\cite{Kn} and \cite{OV} for more details and proofs) and discuss
some aspects of their geometry.

Let ${\mathfrak g}$ be a real semisimple Lie algebra and
${\mathfrak g} = {\mathfrak k} \oplus {\mathfrak p}$ be a
Cartan decomposition of ${\mathfrak g}$.
Let $\theta$ be the corresponding Cartan involution
on ${\mathfrak g}$ and $B$ the Cartan-Killing form on ${\mathfrak g}$.
Then $\langle X,Y \rangle = -B(X,\theta Y)$ is a positive definite inner product
on ${\mathfrak g}$. If $V,W$ are linear subspaces of ${\mathfrak g}$
and $V \subset W$, we denote by $W \ominus V$ the orthogonal complement
of $V$ in $W$ with respect to the inner product,
that is, $W \ominus V = \{ w \in W \mid \langle w,v \rangle = 0\
{\rm\ for\ all\ }v \in V\}$.

Let ${\mathfrak a}$ be
a maximal abelian subspace of ${\mathfrak p}$ and denote by ${\mathfrak a}^\ast$
the dual space of ${\mathfrak a}$. For each $\alpha \in {\mathfrak a}^\ast$
we define ${\mathfrak g}_{\alpha} = \{X \in {\mathfrak g} \mid
[H,X] = \alpha(H)X\ {\rm for\ all\ }H \in {\mathfrak a}\}$.
If $\alpha \neq 0$ and ${\mathfrak g}_\alpha \neq \{0\}$, then $\alpha$
is a restricted root and ${\mathfrak g}_\alpha$ a restricted root
space of ${\mathfrak g}$ with respect to ${\mathfrak a}$. We denote by
$\Sigma$ the set of restricted roots with respect to ${\mathfrak a}$.
The subspace ${\mathfrak g}_0$ coincides with ${\mathfrak k}_0 \oplus
{\mathfrak a}$, where ${\mathfrak k}_0$ is the centralizer of ${\mathfrak a}$
in ${\mathfrak k}$. We recall that ${\mathfrak k}_0 = \{0\}$ if and only if
${\mathfrak g}$ is a split real form of its complexification
${\mathfrak g}^{\mathbb C}$. The direct sum decomposition
$$ {\mathfrak g} = {\mathfrak g}_0 \oplus
\left(\bigoplus_{\alpha \in \Sigma} {\mathfrak g}_{\alpha}\right)
$$ is the restricted root space decomposition of ${\mathfrak g}$
with respect to ${\mathfrak a}$. For each $\alpha \in \Sigma$ we define
the root vector $H_\alpha \in {\mathfrak a}$ corresponding to $\alpha$
by the equation $\alpha(H) = \langle H_\alpha , H \rangle$
for all $H \in {\mathfrak a}$.

Let $\{\alpha_1,\ldots,\alpha_r\} = \Lambda \subset \Sigma$ be a set of simple roots of $\Sigma$,
and denote by $\Sigma^+$ the corresponding set of all positive roots
in $\Sigma$. The subalgebra
$$
{\mathfrak n} = \bigoplus_{\alpha \in
\Sigma^+} {\mathfrak g}_{\alpha}
$$
is nilpotent and
${\mathfrak g} = {\mathfrak k} \oplus {\mathfrak a} \oplus {\mathfrak n}$
is an Iwasawa decomposition of ${\mathfrak g}$.

We will now associate to each subset $\Phi$
of $\Lambda$ a parabolic subalgebra ${\mathfrak q}_\Phi$ of ${\mathfrak g}$.
Let $\Phi$ be a subset of $\Lambda$.
We denote by $\Sigma_\Phi$ the
root subsystem of $\Sigma$ generated by $\Phi$, that is,
$\Sigma_\Phi$ is the intersection of $\Sigma$ and the linear span
of $\Phi$, and put $\Sigma_\Phi^+ = \Sigma_\Phi \cap \Sigma^+$.
We define a reductive subalgebra ${\mathfrak l}_\Phi$ of ${\mathfrak g}$ and a nilpotent subalgebra
${\mathfrak n}_\Phi$ of ${\mathfrak g}$ by
$$
{\mathfrak l}_\Phi = {\mathfrak g}_0 \oplus
\left(\bigoplus_{\alpha \in \Sigma_\Phi} {\mathfrak
g}_{\alpha}\right)
\ \ {\rm and}\ \
{\mathfrak n}_\Phi = \bigoplus_{\alpha \in
\Sigma^+\setminus \Sigma_\Phi^+} {\mathfrak g}_{\alpha}.
$$
Let
$$
{\mathfrak a}_\Phi = \bigcap_{\alpha \in \Phi} {\rm ker}\,\alpha
$$
be the split component of ${\mathfrak l}_\Phi$ and define
${\mathfrak a}^\Phi = {\mathfrak a} \ominus {\mathfrak a}_\Phi$.
Then ${\mathfrak a}_\Phi$ is an abelian subalgebra of ${\mathfrak g}$ and ${\mathfrak l}_\Phi$
is the centralizer and the normalizer of ${\mathfrak a}_\Phi$ in ${\mathfrak g}$.
Since $[{\mathfrak l}_\Phi,{\mathfrak n}_\Phi] \subset {\mathfrak n}_\Phi$,
$$
{\mathfrak q}_\Phi = {\mathfrak l}_\Phi \oplus {\mathfrak n}_\Phi
$$
is a subalgebra of ${\mathfrak g}$,
the so-called parabolic subalgebra of ${\mathfrak g}$
associated with the subsystem $\Phi$ of $\Lambda$.
The decomposition ${\mathfrak q}_\Phi = {\mathfrak l}_\Phi \oplus {\mathfrak n}_\Phi$
is the Chevalley decomposition of the parabolic subalgebra ${\mathfrak q}_\Phi$.

We now define a reductive subalgebra ${\mathfrak m}_\Phi$ of ${\mathfrak g}$ by
${\mathfrak m}_\Phi = {\mathfrak l}_\Phi \ominus {\mathfrak a}_\Phi$.
The subalgebra ${\mathfrak m}_\Phi$ normalizes ${\mathfrak a}_\Phi \oplus {\mathfrak n}_\Phi$, and
${\mathfrak g}_\Phi = [{\mathfrak m}_\Phi,{\mathfrak m}_\Phi] =
[{\mathfrak l}_\Phi,{\mathfrak l}_\Phi]$
is a semisimple subalgebra of ${\mathfrak g}$. The center ${\mathfrak z}_\Phi$ of ${\mathfrak m}_\Phi$
is contained in ${\mathfrak k}_0$ and induces the direct sum decomposition
${\mathfrak m}_\Phi = {\mathfrak z}_\Phi \oplus {\mathfrak g}_\Phi$.
The decomposition
$$
{\mathfrak q}_\Phi = {\mathfrak m}_\Phi \oplus {\mathfrak a}_\Phi \oplus {\mathfrak n}_\Phi
$$
is the Langlands decomposition of the parabolic subalgebra ${\mathfrak q}_\Phi$.

For $\Phi = \emptyset$ we have ${\mathfrak l}_\emptyset =
{\mathfrak g}_0$, ${\mathfrak m}_\emptyset = {\mathfrak k}_0$,
${\mathfrak a}_\emptyset = {\mathfrak a}$ and ${\mathfrak n}_\emptyset =
{\mathfrak n}$. In this case ${\mathfrak q}_\emptyset = {\mathfrak
k}_0 \oplus {\mathfrak a} \oplus {\mathfrak n} =
{\mathfrak g}_0 \oplus {\mathfrak n}$ is a minimal
parabolic subalgebra of ${\mathfrak g}$. For $\Phi = \Lambda$ we
obtain ${\mathfrak l}_\Lambda = {\mathfrak m}_\Lambda = {\mathfrak g}$ and
${\mathfrak a}_\Lambda = {\mathfrak n}_\Lambda = \{0\}$. The proper maximal parabolic
subalgebras of ${\mathfrak g}$ are precisely those parabolic subalgebras for which
the cardinality $|\Lambda \setminus \Phi|$ of $\Lambda \setminus \Phi$
is equal to one. The proper maximal parabolic subalgebras
can therefore be parametrized by the simple roots in $\Lambda$.

Each parabolic subalgebra of ${\mathfrak g}$ is conjugate in ${\mathfrak g}$ to
${\mathfrak q}_\Phi$ for some subset $\Phi$ of $\Lambda$.
The set of conjugacy classes of parabolic subalgebras of ${\mathfrak g}$ therefore
has $2^r$ elements, where $r = |\Lambda|$ is the real rank of ${\mathfrak g}$.
Two parabolic subalgebras ${\mathfrak q}_{\Phi_1}$ and
${\mathfrak q}_{\Phi_2}$ of ${\mathfrak g}$ are conjugate in the full
automorphism group ${\rm Aut}({\mathfrak g})$ of ${\mathfrak g}$
if and only if there exists an automorphism $F$ of the Dynkin diagram
associated to $\Lambda$ with $F(\Phi_1) = \Phi_2$.

For each $\alpha \in \Sigma$
we define
${\mathfrak k}_\alpha = {\mathfrak k} \cap ({\mathfrak g}_{-\alpha} \oplus {\mathfrak g}_\alpha)$ and
${\mathfrak p}_\alpha = {\mathfrak p} \cap ({\mathfrak g}_{-\alpha} \oplus {\mathfrak g}_\alpha)$.
Then we have ${\mathfrak k}_{-\alpha} = {\mathfrak k}_\alpha$, ${\mathfrak p}_{-\alpha} = {\mathfrak p}_\alpha$
and ${\mathfrak k}_\alpha \oplus {\mathfrak p}_\alpha = {\mathfrak g}_{-\alpha} \oplus {\mathfrak g}_\alpha$
for all $\alpha \in \Sigma$.
It is easy to see that the subspaces
$$
{\mathfrak p}_\Phi = {\mathfrak l}_\Phi \cap {\mathfrak p} =
{\mathfrak a} \oplus \left( \bigoplus_{\alpha \in \Sigma_{\Phi}}
{\mathfrak p}_\alpha \right)
\ {\rm and}\
{\mathfrak b}_\Phi = {\mathfrak m}_\Phi \cap {\mathfrak p} =
{\mathfrak g}_\Phi \cap {\mathfrak p} =
{\mathfrak a}^\Phi \oplus \left( \bigoplus_{\alpha \in \Sigma_{\Phi}}
{\mathfrak p}_\alpha \right)
$$
are Lie triple systems in ${\mathfrak p}$.
We define a subalgebra ${\mathfrak k}_\Phi$ of ${\mathfrak k}$ by
$$
{\mathfrak k}_\Phi = {\mathfrak q}_\Phi \cap {\mathfrak k} = {\mathfrak l}_\Phi \cap {\mathfrak k}
= {\mathfrak m}_\Phi \cap {\mathfrak k} = {\mathfrak k}_0 \oplus \left( \bigoplus_{\alpha \in \Sigma_{\Phi}}
{\mathfrak k}_\alpha \right).
$$
Then we have
$$
[{\mathfrak k}_\Phi , {\mathfrak m}_\Phi ] \subset {\mathfrak m}_\Phi\ ,\
[{\mathfrak k}_\Phi , {\mathfrak a}_\Phi ] = \{0\}\ ,\
[{\mathfrak k}_\Phi , {\mathfrak n}_\Phi ] \subset {\mathfrak n}_\Phi\ .
$$
These three relations will be important for our understanding of cohomogeneity one actions on $M$.
Moreover, ${\mathfrak g}_\Phi = ({\mathfrak g}_\Phi \cap {\mathfrak k}_\Phi) \oplus {\mathfrak b}_\Phi$
is a Cartan decomposition of the semisimple subalgebra ${\mathfrak g}_\Phi$
of ${\mathfrak g}$ and ${\mathfrak a}^\Phi$ is a maximal abelian subspace
of ${\mathfrak b}_\Phi$. If we define
$({\mathfrak g}_\Phi)_0 = ({\mathfrak g}_\Phi \cap {\mathfrak k}_0) \oplus {\mathfrak a}^\Phi$,
then
${\mathfrak g}_\Phi = ({\mathfrak g}_\Phi)_0 \oplus
\left(\bigoplus_{\alpha \in \Sigma_\Phi} {\mathfrak
g}_{\alpha}\right)$
is the restricted root space decomposition of ${\mathfrak g}_\Phi$ with respect to ${\mathfrak a}^\Phi$
and $\Phi$ is the corresponding set of simple roots. Since ${\mathfrak m}_\Phi = {\mathfrak z}_\Phi
\oplus {\mathfrak g}_\Phi$ and ${\mathfrak z}_\Phi \subset {\mathfrak k}_0$, we see that
${\mathfrak g}_\Phi \cap {\mathfrak k}_0 = {\mathfrak k}_0 \ominus {\mathfrak z}_\Phi$.

We now relate these algebraic constructions to the geometry of symmetric spaces
of noncompact type.
Let $M = G/K$ be the connected Riemannian symmetric space of noncompact type
associated with the pair $({\mathfrak g},{\mathfrak k})$.
The Riemannian metric on $M$ is the one which is induced from the ${\rm Ad}(K)$-invariant inner product
$\langle \cdot , \cdot \rangle$ on $\mathfrak p$. Then
$G = I^o(M)$ is the connected component of the isometry group of $M$
containing the identity and $K$ is a maximal compact subgroup of $G$.
The Lie algebra of $G$ and $K$ coincides with ${\mathfrak g}$ and ${\mathfrak k}$,
respectively. We denote by $o \in M$ the unique fixed point of $K$, that is,
$o$ is the point in $M$ for which the stabilizer of $G$ at $o$ coincides with $K$.
We identify the subspace ${\mathfrak p}$ in the Cartan decomposition
${\mathfrak g} = {\mathfrak k} \oplus {\mathfrak p}$ with the tangent space
$T_oM$ of $M$ at $o$ in the usual way. The rank of the symmetric space
$M$ coincides with $r = |\Lambda|$.

Let ${\rm Exp} : {\mathfrak g} \to G$ be the Lie exponential map of ${\mathfrak g}$.
Then $A = {\rm Exp}({\mathfrak a})$ and $N = {\rm Exp}({\mathfrak n})$ is a simply
connected closed subgroup of $G$ with Lie algebra ${\mathfrak a}$ and ${\mathfrak n}$,
respectively, $A$ is abelian and $N$ is nilpotent.
The orbit $A \cdot o$ is an $r$-dimensional Euclidean space ${\mathbb E}^r$
embedded totally geodesically into $M$, and the orbit $N\cdot o$ is a horocycle
in $M$. The Iwasawa decomposition
${\mathfrak g} = {\mathfrak k} \oplus {\mathfrak a} \oplus {\mathfrak n}$ of ${\mathfrak g}$
induces an Iwasawa decomposition $G = KAN$ of $G$. The solvable Lie group $AN$ acts
simply transitively on the symmetric space $M$.

Let $\Phi$ be a subset of $\Lambda$ and $r_\Phi = |\Phi|$. We denote by $A_\Phi$
the connected abelian subgroup of $G$ with Lie algebra ${\mathfrak a}_\Phi$ and by
$N_\Phi$ the connected nilpotent subgroup of $G$ with Lie algebra ${\mathfrak n}_\Phi$.
The centralizer $L_\Phi = Z_G({\mathfrak a}_\Phi)$ of ${\mathfrak a}_\Phi$ in $G$ is a reductive subgroup
of $G$ with Lie algebra ${\mathfrak l}_\Phi$.
Moreover, $L_\Phi$ normalizes $N_\Phi$, and hence $Q_\Phi = L_\Phi N_\Phi$ is a subgroup of
$G$ with Lie algebra ${\mathfrak q}_\Phi$. The subgroup $Q_\Phi$ coincides with the normalizer
$N_G({\mathfrak l}_\Phi \oplus {\mathfrak n}_\Phi)$ of ${\mathfrak l}_\Phi \oplus {\mathfrak n}_\Phi$
in $G$, and hence $Q_\Phi$ is a closed subgroup of $G$.
The subgroup $Q_\Phi$ is the parabolic subgroup of $G$ associated with the subsystem $\Phi$
of $\Lambda$. We denote by $Q_\Phi^o$ the connected component of $Q_\Phi$ containing
the identity transformation.

Let $G_\Phi$ be the connected subgroup of $G$ with Lie algebra
${\mathfrak g}_\Phi$. Since ${\mathfrak g}_\Phi$ is semisimple, $G_\Phi$ is a semisimple subgroup of $G$.
The intersection $K_\Phi = L_\Phi \cap K$
is a maximal compact subgroup of $L_\Phi$ and ${\mathfrak k}_\Phi$ is the Lie algebra of $K_\Phi$.
The adjoint group ${\rm Ad}(L_\Phi)$ normalizes ${\mathfrak g}_\Phi$, and consequently
$M_\Phi = K_\Phi G_\Phi$ is a subgroup of $L_\Phi$. One can show that
$M_\Phi$ is a closed reductive subgroup of $L_\Phi$,
$K_\Phi$ is a maximal compact subgroup of $M_\Phi$,
and the center $Z_\Phi$ of $M_\Phi$ is a compact subgroup of $K_\Phi$. The Lie algebra
of $M_\Phi$ is ${\mathfrak m}_\Phi$ and $L_\Phi$ is isomorphic to the
Lie group direct product $M_\Phi \times A_\Phi$. The multiplication
$M_\Phi \times A_\Phi \times N_\Phi \to Q_\Phi$
is an analytic diffeomorphism, and the group structure is given by
$$
(m,a,n)(m^\prime,a^\prime,n^\prime) = (mm^\prime,aa^\prime,(m^\prime a^\prime)^{-1}n(m^\prime a^\prime) n^\prime).
$$
The parabolic subgroup $Q_\Phi$ acts transitively on $M$ and the isotropy subgroup at $o$
is $K_\Phi$, that is, $M = Q_\Phi/K_\Phi$.

Since ${\mathfrak g}_\Phi = ({\mathfrak g}_\Phi \cap {\mathfrak k}_\Phi) \oplus {\mathfrak b}_\Phi$
is a Cartan decomposition of the semisimple subalgebra ${\mathfrak g}_\Phi$, we have
$[{\mathfrak b}_\Phi,{\mathfrak b}_\Phi] = {\mathfrak g}_\Phi \cap {\mathfrak k}_\Phi$.
Thus $G_\Phi$ is the connected closed subgroup of $G$ with Lie algebra
$[{\mathfrak b}_\Phi,{\mathfrak b}_\Phi] \oplus {\mathfrak b}_\Phi$.
Since ${\mathfrak b}_\Phi$ is a Lie triple system in ${\mathfrak p}$, the orbit
$B_\Phi = G_\Phi \cdot o$ of the $G_\Phi$-action on $M$ containing $o$
is a connected totally geodesic submanifold of $M$ with $T_oB_\Phi = {\mathfrak b}_\Phi$.
If $\Phi = \emptyset$, then $B_\Phi = \{o\}$, otherwise
$B_\Phi$ is a Riemannian symmetric space of noncompact type and ${\rm rank}(B_\Phi) = r_\Phi$, and
$$B_\Phi = G_\Phi \cdot o = G_\Phi/(G_\Phi\cap K_\Phi) = M_\Phi \cdot o = M_\Phi/K_\Phi.$$
The submanifold $B_\Phi$ is also known as a boundary component of $M$ in
the context of the maximal Satake compactification of $M$ (see e.g.\ \cite{BJ}).

Clearly, ${\mathfrak a}_\Phi$ is a Lie triple system as well, and the corresponding
totally geodesic submanifold is a Euclidean space
$${\mathbb E}^{r-r_\Phi} = A_\Phi \cdot o.$$

Finally, ${\mathfrak p}_\Phi = {\mathfrak b}_\Phi \oplus {\mathfrak a}_\Phi$ is a Lie triple system, and the
corresponding totally geodesic submanifold $F_\Phi$ is the symmetric space
$$
F_\Phi = L_\Phi \cdot o = L_\Phi/K_\Phi = (M_\Phi \times A_\Phi)/K_\Phi =
 B_\Phi \times {\mathbb E}^{r-r_\Phi}.
$$

The analytic diffeomorphism $M_\Phi \times A_\Phi \times N_\Phi \to Q_\Phi$ induces an
analytic diffeomorphism
$$
B_\Phi \times A_\Phi \times N_\Phi \to M, (m \cdot o,a,n) \mapsto (man)\cdot o,
$$
known as a horospherical decomposition
of the symmetric space $M$. The action of $Q_\Phi$ on $M$ is given by
$$
Q_\Phi \times M \to M , ((m,a,n),(m^\prime \cdot o,a^\prime,n^\prime)) \mapsto ((mm^\prime) \cdot o,aa^\prime,
(m^\prime a^\prime)^{-1}n(m^\prime a^\prime) n^\prime).
$$

\section{Maximal reductive and parabolic subgroups}

In this section we relate cohomogeneity one actions on $M$ to actions of reductive and parabolic subgroups of $G$.

\begin{prop} \label{redtotgeod}
Every connected proper maximal reductive subgroup $L$ of $G$ has a totally
geodesic orbit $W$ in $M$ with $\dim W < \dim M$. In particular,
$L$ cannot act transitively on $M$.
\end{prop}

\begin{proof}
Let ${\mathfrak l}$ be the Lie algebra of $L$. As ${\mathfrak g}$
is algebraic (see e.g.\ \cite{OV}, p.~29, Corollary~4) and ${\mathfrak l}$
is maximal in ${\mathfrak g}$, ${\mathfrak l}$ is an algebraic
subalgebra of ${\mathfrak g}$. Since
${\mathfrak l}$ is a reductive algebraic subalgebra
of ${\mathfrak g}$, there exists a Cartan decomposition
${\mathfrak g} = {\mathfrak k} \oplus {\mathfrak p}$ of
${\mathfrak g}$ such that ${\mathfrak l} = ({\mathfrak l} \cap
{\mathfrak k}) \oplus ({\mathfrak l} \cap {\mathfrak p})$ (see
e.g.\ \cite{OV}, p.~207, Theorem~3.6). Let $K$ be the maximal
compact subgroup of $G$ with Lie algebra ${\mathfrak k}$, and let
$o \in M$ be the fixed point of $K$. Then the orbit $W = L \cdot
o$ of $L$ through $o$ is a totally geodesic submanifold of $M$
(see e.g.\ \cite{BCO}, Proposition~9.1.2). Assume that $\dim W =
\dim M$, which means that
${\mathfrak p} = {\mathfrak l} \cap {\mathfrak p} \subset {\mathfrak l}$.
Since ${\mathfrak g}$ is semisimple and contains no nonzero compact ideals,
we have $[{\mathfrak p},{\mathfrak p}] = {\mathfrak k}$
(see e.g.\ \cite{OV}, p.~145, Proposition~3.5).
This implies ${\mathfrak k} = [{\mathfrak
p},{\mathfrak p}] \subset [{\mathfrak l},{\mathfrak l}] \subset
{\mathfrak l}$. Altogether this gives ${\mathfrak g} = {\mathfrak k} \oplus
{\mathfrak p} \subset {\mathfrak l}$ and hence ${\mathfrak l} =
{\mathfrak g}$. As $G$ is connected, this contradicts the
assumption that $L$ is a proper subgroup of $G$, and we conclude
$\dim W < \dim M$.
\end{proof}

{\it Remarks.} 1. It was shown by Karpelevic \cite{Kar}
that every connected semisimple subgroup of $G$ has a totally geodesic orbit
in $M$. This follows also from Theorem 6 proved by Mostow in \cite{Mo1}.
A geometric proof for the semisimple case was recently given by
Di Scala and Olmos in \cite{DO}.

\medskip
2. The corresponding statement for Riemannian symmetric
spaces of compact type is not true. Consider for example the
$4$-dimensional sphere $S^4 = SO(5)/SO(4)$ as a subset of the
$5$-dimensional real vector space of all symmetric $(3 \times
3)$-matrices with real coefficients and trace zero. By considering
the action of $SO(3)$ on such matrices by conjugation one gets a
cohomogeneity one action on $S^4$ with no totally geodesic orbit.
The two singular orbits of this action are congruent to the
Veronese embedding of the real projective plane ${\mathbb R}P^2$
into $S^4$. In the compact case there also exist connected proper
reductive subgroups which act transitively. For example, $SU(n)$
($n \geq 2$) is a connected proper reductive subgroup of $SO(2n)$
which acts transitively on $SO(2n)/SO(2n-1) = S^{2n-1}$.

\begin{thm} \label{maximal}
Let $M$ be a connected Riemannian symmetric space of noncompact
type and $H$ a connected subgroup of $G = I^o(M)$ acting on $M$
with cohomogeneity one. Let $L$ be a connected proper maximal subgroup
of $G$ with $H \subset L$. Then one of the following statements
holds:
\begin{itemize}
\item[(1)] $L$ is a reductive subgroup of $G$, the actions of $H$ and $L$ are
orbit equivalent, and the action of
$H$ on $M$ has a totally geodesic orbit $W$. Moreover, if $M$ is irreducible and
$M \neq {\mathbb R}H^n = SO^o_{1,n}/SO_n$, then $W$ is a singular orbit.
\item[(2)] $L$ is the identity component of a parabolic
subgroup of $G$.
\end{itemize}
\end{thm}

\begin{proof} We denote by
${\mathfrak l}$ the Lie algebra of $L$, by ${\mathfrak r}$ the
radical of ${\mathfrak l}$, and by ${\mathfrak n}$ the nilradical
of ${\mathfrak l}$. It is a well-known consequence of Lie's
Theorem on solvable Lie algebras that $[{\mathfrak l},{\mathfrak
r}] \subset {\mathfrak n}$ (see e.g.\ \cite{Kn}, Corollary 1.41).
Mostow has shown (see proof of Theorem 3.1 in \cite{Mo2}) that
the nilradical ${\mathfrak n}$ is trivial if and only if $L$ is unimodular.

Let us first assume that $L$ is unimodular. Then $[{\mathfrak
l},{\mathfrak r}] = 0$, which implies that ${\mathfrak r}$ is
contained in the center ${\mathfrak z}$ of ${\mathfrak l}$. As the
center of a Lie algebra is always contained in the radical of the
Lie algebra, we conclude that the radical ${\mathfrak r}$ of
${\mathfrak l}$ coincides with the center ${\mathfrak z}$ of
${\mathfrak l}$. Therefore ${\mathfrak l}$ is a reductive Lie
algebra. As $H \subset L$, the orbits of the action of $H$ are
contained in the orbits of the action of $L$. However, $L$ cannot
act transitively on $M$ (see Proposition \ref{redtotgeod}) and hence
must act on $M$ with cohomogeneity one. Since
both $L$ and $H$ are connected, the orbits of $H$ and
$L$ must therefore coincide, and Proposition \ref{redtotgeod} implies that
$H$ has a totally geodesic orbit $W$ with $\dim W < \dim M$.

The real hyperbolic spaces ${\mathbb R}H^n$, $n \geq 2$, are the
only irreducible Riemannian symmetric spaces of noncompact type
which have a totally geodesic hypersurface (see e.g.\ \cite{Iw}).
Therefore, if $M$ is irreducible and
$M \neq {\mathbb R}H^n = SO^o_{1,n}/SO_n$, the
totally geodesic orbit $W$ must be a singular orbit of the action.

If $L$ is not unimodular, then ${\mathfrak l}$ is a parabolic
subalgebra of ${\mathfrak g}$ by a result of
Mostow (\cite{Mo2}), and hence $L$ is the identity component
of a parabolic subgroup of $G$. Therefore $L$ is conjugate to $Q_\Phi^o$ for some
subset $\Phi$ of $\Lambda$, and since $L$ is a maximal proper subgroup of $G$, we have $|\Phi| = r-1$.
\end{proof}

{\it Remark.} The maximal reductive nonsemisimple subalgebras of
real semisimple Lie algebras have been classified by Tao \cite{Tao}.

\smallskip
In view of Theorem \ref{maximal} we now  consider more thoroughly
 the case when ${\mathfrak l}$ is a parabolic subalgebra of ${\mathfrak g}$.

\section{The parabolic case}
\label{parabolic case}

In this section we assume that ${\mathfrak h}$ is contained in a parabolic
subalgebra ${\mathfrak l}$ of ${\mathfrak g}$. From Section \ref{parabolic} we know
that ${\mathfrak l}$ is conjugate to ${\mathfrak q}_\Phi$ for some subset $\Phi$
of $\Lambda$. Without loss of generality we assume that ${\mathfrak l} = {\mathfrak q}_\Phi$.
Now consider the Langlands decomposition
$$
{\mathfrak q}_\Phi = {\mathfrak m}_\Phi \oplus {\mathfrak a}_\Phi \oplus {\mathfrak n}_\Phi
$$
of ${\mathfrak q}_\Phi$
and the corresponding horospherical decomposition
$$ M \cong B_\Phi \times {\mathbb E}^{r-r_\Phi} \times N_\Phi \cong
B_\Phi \times A_\Phi \times N_\Phi$$
of $M$. Note that for the second congruence we identify the Euclidean space ${\mathbb E}^{r-r_\Phi}$
and the abelian Lie group $A_\Phi$ via the simple transitive action of $A_\Phi$ on ${\mathbb E}^{r-r_\Phi}$.
We now construct two types of cohomogeneity one actions from
the Langlands or horospherical decomposition.

\subsection{Canonical extensions from boundary components}

Let $H_\Phi$ be a connected subgroup of $I(B_\Phi)$
and denote by ${\mathfrak h}_\Phi$ the Lie algebra
of $H_\Phi$. Since ${\mathfrak h}_\Phi \subset {\mathfrak g}_\Phi \subset {\mathfrak m}_\Phi$
and ${\mathfrak m}_\Phi$ normalizes ${\mathfrak a}_\Phi \oplus {\mathfrak n}_\Phi$,
we see that
$$
{\mathfrak h}_\Phi^\Lambda = {\mathfrak h}_\Phi \oplus {\mathfrak a}_\Phi \oplus {\mathfrak n}_\Phi
$$
is a subalgebra of ${\mathfrak m}_\Phi \oplus {\mathfrak a}_\Phi \oplus {\mathfrak n}_\Phi
= {\mathfrak q}_\Phi \subset {\mathfrak g}$.
We call the connected subgroup $H_\Phi^\Lambda$ of
the parabolic subgroup $Q_\Phi$ of $G$ with Lie algebra
${\mathfrak h}_\Phi^\Lambda$ the {\it canonical extension} of $H_\Phi$ from
the boundary component $B_\Phi$ to the symmetric space $M$.

By construction we have
${\mathfrak h}_\Phi^\Lambda \cap {\mathfrak k} = {\mathfrak h}_\Phi \cap {\mathfrak k}$, and
the normal space at $o$ in $T_oM$ of the orbit $H_\Phi^\Lambda \cdot o$ coincides
with the normal space at $o$ in $T_oB_\Phi$ of the orbit $H_\Phi \cdot o$. This implies that
the slice representations of $H_\Phi^\Lambda$ and $H_\Phi$ at $o$ coincide. Therefore,
the cohomogeneity of the action of $H_\Phi^\Lambda$ on $M$ coincides with
the cohomogeneity of the action of $H_\Phi$ on $B_\Phi$. We therefore conclude:

\begin{prop}
\label{subdiagram} Let $M$ be a connected Riemannian symmetric space of noncompact type
and let $B_\Phi$ be a boundary component of $M$.
Then every cohomogeneity one action on $B_\Phi$
has a canonical extension to a cohomogeneity one action on $M$.
\end{prop}

{\it Remark.} If $\Phi \subset \Psi \subset \Lambda$ are nonempty proper subsets of each other,
then clearly $B_\Phi \subset B_\Psi \subset M$ are proper totally geodesic
submanifolds of each other, and $B_\Phi$ is a boundary component
of the semisimple symmetric space $B_\Psi$. Let $H_\Phi$ be a
connected subgroup of $I(B_\Phi)$ acting
on $B_\Phi$ with cohomogeneity one, and denote by $H_\Phi^\Psi$ the connected subgroup of $I(B_\Psi)$
obtained by canonical extension of the $H_\Phi$-action from $B_\Phi$ to $B_\Psi$,
by $H_\Phi^\Lambda$ the connected subgroup of $I(M)$ obtained by
canonical extension of the $H_\Phi$-action from $B_\Phi$ to $M$, and by
$H_\Psi^\Lambda$ the connected subgroup of $I(M)$ obtained by
canonical extension of the $H_\Phi^\Psi$-action from $B_\Psi$ to $M$.
It follows from the construction that the corresponding Lie algebras satisfy
${\mathfrak h}_\Phi^\Lambda = {\mathfrak h}_\Psi^\Lambda$, and therefore $H_\Phi^\Lambda = H_\Psi^\Lambda$
as both groups are connected. This shows that for the classification
of cohomogeneity one actions on $M$ obtained by canonical extensions one can restrict
to canonical extensions of cohomogeneity one actions on boundary components
$B_\Phi$ of rank $r-1$, that is, those boundary components $B_\Phi$ obtained from subsets $\Phi$ of
$\Lambda$ with $|\Phi| = r-1$.

\medskip
We now investigate in how far canonical extensions preserve orbit equivalence of cohomogeneity one actions.

\begin{prop} \label{orbit equivalence}
Let $M$ be a connected Riemannian symmetric space of noncompact type
and let $B_\Phi$ be a boundary component of $M$. Let $H^1_\Phi,H^2_\Phi$ be two connected closed subgroups of $I(B_\Phi)$ which act on $B_\Phi$ with cohomogeneity one. Assume that these two actions are orbit equivalent by an isometry $f \in I^o(B_\Phi)$. Then the two cohomogeneity one actions on $M$ which are obtained by canonical extension of $H^1_\Phi$ and $H^2_\Phi$ are orbit equivalent.
\end{prop}

\begin{proof}
Since $G_\Phi$ is a connected semisimple Lie group acting transitively on $B_\Phi$, we must have $I^o(B_\Phi) \subset G_\Phi$. Since $G_\Phi \subset M_\Phi \subset L_\Phi \subset Q_\Phi \subset G$, the isometry $f$ extends canonically to an isometry $F$ in the parabolic subgroup $Q_\Phi$ of $G$. The horospherical decomposition $M = B_\Phi \times A_\Phi \times N_\Phi$ shows that
$$
F((H_\Phi^1)^\Lambda \cdot p) = f(H_\Phi^1 \cdot p) \times A_\Phi \times N_\Phi.
$$
By assumption, we have
$$
f(H_\Phi^1 \cdot p) \times A_\Phi \times N_\Phi = (H_\Phi^2 \cdot f(p)) \times A_\Phi \times N_\Phi.
$$
Recall that we have an analytic diffeomorphism $M_\Phi \times A_\Phi \times N_\Phi \to Q_\Phi$,
and accordingly we write $F = (\bar{m},1,1)$ with $\bar{m} \in G_\Phi$.
Since the group structure is given by
$$
(m,a,n)(m^\prime,a^\prime,n^\prime) = (mm^\prime,aa^\prime,(m^\prime a^\prime)^{-1}n(m^\prime a^\prime) n^\prime),
$$
we obtain
\begin{eqnarray*}
(H_\Phi^2)^\Lambda \cdot F(p) & = & \{ (m,a,n)(\bar{m},1,1) \cdot p \mid m \in H_\Phi^2,\ a \in A_\Phi,\ n \in N_\Phi\}\\
& = & \{ (m\bar{m},a,\bar{m}^{-1}n\bar{m}) \cdot p \mid m \in H_\Phi^2,\ a \in A_\Phi,\ n \in N_\Phi\}\\
& = & \{ (m\bar{m},a,n) \cdot p \mid m \in H_\Phi^2,\ a \in A_\Phi,\ n \in N_\Phi\} \\
& = & (H_\Phi^2 \cdot f(p)) \times A_\Phi \times N_\Phi.
\end{eqnarray*}
Altogether we see that $F((H_\Phi^1)^\Lambda \cdot p) = (H_\Phi^2)^\Lambda \cdot F(p)$, which means that the two cohomogeneity one actions on $M$ obtained by canonical extension of $H^1_\Phi$ and $H^2_\Phi$ are orbit equivalent.
\end{proof}

\medskip
The following example shows that we cannot weaken the assumption in Proposition \ref{orbit equivalence} from $f \in I^o(B_\Phi)$ to $f \in I(B_\Phi)$.

\medskip
{\it Example.} We consider the symmetric space $M = SL_4({\mathbb R})/SO_4$.
This symmetric space has rank $3$, dimension $9$, and the restricted root system is of type $(A_3)$
with all multiplicities equal to one. We choose $\Phi = \{\alpha_1,\alpha_2\} \subset \Lambda = \{\alpha_1,\alpha_2,\alpha_3\}$. Then we have
$$
{\mathfrak l}_\Phi = {\mathfrak g}_0 \oplus {\mathfrak g}_{\pm \alpha_1} \oplus {\mathfrak g}_{\pm \alpha_2} \oplus {\mathfrak g}_{\pm (\alpha_1 + \alpha_2)} \cong {\mathfrak s}{\mathfrak l}_3({\mathbb R}) \oplus {\mathbb R},
$$
and therefore
$$
F_\Phi = L_\Phi \cdot o \cong SL_3({\mathbb R})/SO_3 \times {\mathbb E}.
$$
The corresponding boundary component $B_\Phi$ is isometric to $SL_3({\mathbb R})/SO_3$. We now define two
subalgebras ${\mathfrak h}_\Phi^1$ and ${\mathfrak h}_\Phi^2$ of ${\mathfrak s}{\mathfrak l}_3({\mathbb R})$ by
$$
{\mathfrak h}^1_\Phi = {\mathfrak a}^\Phi \oplus {\mathfrak g}_{\alpha_2} \oplus {\mathfrak g}_{\alpha_1 + \alpha_2}\ \ {\rm and}\ \
{\mathfrak h}^2_\Phi = {\mathfrak a}^\Phi \oplus {\mathfrak g}_{\alpha_1} \oplus {\mathfrak g}_{\alpha_1 + \alpha_2}.
$$
The corresponding connected subgroups $H^1_\Phi$ and $H^2_\Phi$ of $SL_3({\mathbb R})$ act on the boundary component $B_\Phi$ with cohomogeneity one, and the orbits form a foliation on $B_\Phi$. These two actions are orbit equivalent, and the corresponding isometry is induced by the Dynkin diagram symmetry of $(A_2)$, the restricted root system of $B_\Phi$. We now consider the canonical extensions of these two actions, which are defined by
$$ ({\mathfrak h}_\Phi^1)^\Lambda = {\mathfrak h}_\Phi^1 \oplus {\mathfrak a}_\Phi \oplus {\mathfrak n}_\Phi\ \ {\rm and}\ \ ({\mathfrak h}_\Phi^2)^\Lambda = {\mathfrak h}_\Phi^2 \oplus {\mathfrak a}_\Phi \oplus {\mathfrak n}_\Phi.
 $$
In terms of the root space decomposition of $M$ these two subalgebras are
$$
({\mathfrak h}_\Phi^1)^\Lambda = {\mathfrak a} \oplus ({\mathfrak n} \ominus {\mathfrak g}_{\alpha_1})\ \ {\rm and}\ \ ({\mathfrak h}_\Phi^2)^\Lambda = {\mathfrak a} \oplus ({\mathfrak n} \ominus {\mathfrak g}_{\alpha_2})
$$
The corresponding connected subgroups $(H^1_\Phi)^\Lambda$ and $(H^2_\Phi)^\Lambda$ of
$SL_4({\mathbb R})$ act on the symmetric space $M = SL_4({\mathbb R})/SO_4$ with cohomogeneity one, and the orbits form a foliation on $M$. However, these two actions are not orbit equivalent since there is no corresponding Dynkin diagram symmetry (see \cite{BT1} for details). The reason for this is that the Dynkin diagram symmetry of $(A_2)$ does not extend to a Dynkin diagram symmetry of $(A_3)$.

\subsection{Nilpotent construction} \label{method2}

We now describe our second new method for constructing cohomogeneity one
actions on $M$.

Let $\Phi$ be a subset of $\Lambda$ and
consider the parabolic subalgebra ${\mathfrak q}_\Phi$ and its
Langlands decomposition ${\mathfrak q}_\Phi = {\mathfrak m}_\Phi
\oplus {\mathfrak a}_\Phi \oplus {\mathfrak n}_\Phi$.
Recall that
${\mathfrak q}_\Phi \cap {\mathfrak k} =
{\mathfrak m}_\Phi \cap {\mathfrak k} = {\mathfrak k}_\Phi$,
and we have a canonical isomorphism
$$T_oM \cong {\mathfrak b}_\Phi \oplus {\mathfrak a}_\Phi \oplus {\mathfrak n}_\Phi\ ,\
T_oF_\Phi \cong {\mathfrak b}_\Phi \oplus {\mathfrak a}_\Phi\ ,\
T_oB_\Phi \cong {\mathfrak b}_\Phi.$$
Since $K_\Phi \subset L_\Phi = Z_G({\mathfrak a}_\Phi)$, we have
${\rm Ad}(k)X = X$ for all $k \in K_\Phi$ and $X \in {\mathfrak a}_\Phi$.
Furthermore, since $K_\Phi \subset M_\Phi$ and ${\mathfrak m}_\Phi = {\mathfrak k}_\Phi \oplus
{\mathfrak b}_\Phi$ is a Cartan decomposition, we get
${\rm Ad}(k)({\mathfrak b}_\Phi) = {\mathfrak b}_\Phi$ for all $k \in K_\Phi$.
Eventually, since $K_\Phi \subset L_\Phi$ normalizes $N_\Phi$, we get
${\rm Ad}(k)({\mathfrak n}_\Phi) = {\mathfrak n}_\Phi$ for all $k \in K_\Phi$.
Thus the decomposition $T_oM \cong {\mathfrak b}_\Phi \oplus {\mathfrak a}_\Phi \oplus {\mathfrak n}_\Phi$ is ${\rm Ad}(K_\Phi)$-invariant.

The nilpotent subalgebra ${\mathfrak n}_\Phi$ has a natural gradation which we shall
now describe. Let $H^1, \ldots, H^r \in {\mathfrak a}$ be the dual
vectors of $\alpha_1,\ldots,\alpha_r$, that is, define $H^1,\ldots,H^r \in {\mathfrak a}$ by
$\alpha_\nu(H^\mu) = \delta_{\nu\mu}$. Define
$$
H^\Phi = \sum_{\alpha_i \in \Lambda \setminus \Phi} H^i
$$
and
$m_\Phi = \tilde{\alpha}(H^\Phi)$, where $\tilde{\alpha}$ is the highest root
in $\Sigma^+$. For each $\alpha \in \Sigma^+$ we have $\alpha(H^\Phi) \in \{0,\ldots,m_\Phi\}$,
and $\alpha \in \Sigma_\Phi^+$ if and only if $\alpha(H^\Phi) = 0$.
For each $\nu \in \{1,\ldots,m_\Phi\}$ we define a subspace ${\mathfrak n}_\Phi^\nu$
of ${\mathfrak n}_\Phi$ by
$$
{\mathfrak n}_\Phi^\nu = \bigoplus_{\substack{\alpha \in \Sigma^+ \setminus \Sigma_\Phi^+ \\ \alpha(H^\Phi)=\nu}}
{\mathfrak g}_\alpha .
$$
Then
$$
{\mathfrak n}_\Phi = \bigoplus_{\nu = 1}^{m_\Phi} {\mathfrak n}_\Phi^\nu
$$
is an ${\rm Ad}(K_\Phi)$-invariant gradation of ${\mathfrak g}$.
This gradation is generated by ${\mathfrak n}_\Phi^1$, which means that
$[{\mathfrak n}_\Phi^1,{\mathfrak n}_\Phi^\nu] = {\mathfrak n}_\Phi^{\nu+1}$
holds for all $\nu \in \{1,\ldots,m_\Phi-1\}$ (see \cite{KA}). It is clear
that ${\mathfrak n}_\Phi$ is abelian if and only if $m_\Phi = 1$.

Assume that $\dim {\mathfrak n}_\Phi^1 \geq 2$ and
let ${\mathfrak v}$ be a subspace of ${\mathfrak n}_\Phi^1$ with
$\dim {\mathfrak v} \geq 2$. Since $[{\mathfrak n}_\Phi , {\mathfrak n}_\Phi ]
= {\mathfrak n}_\Phi \ominus {\mathfrak n}_\Phi^1$, we see that
$$
{\mathfrak n}_{\Phi,{\mathfrak v}} = {\mathfrak n}_\Phi \ominus {\mathfrak v}
$$
is a subalgebra of ${\mathfrak n}_\Phi$.
Let $N_{\Phi,{\mathfrak v}}$ be the connected subgroup of $N_\Phi$
with Lie algebra ${\mathfrak n}_{\Phi,{\mathfrak v}}$ and
$N_{L_\Phi}({\mathfrak n}_{\Phi,{\mathfrak v}})$ be the normalizer
of ${\mathfrak n}_{\Phi,{\mathfrak v}}$ in $L_\Phi$. Since $K_\Phi = L_\Phi \cap K$,
the normalizer $N_{K_\Phi}({\mathfrak n}_{\Phi,{\mathfrak v}})$ of
${\mathfrak n}_{\Phi,{\mathfrak v}}$ in $K_\Phi$ coincides with
$N_{L_\Phi}({\mathfrak n}_{\Phi,{\mathfrak v}}) \cap K$. Moreover,
since ${\rm Ad}(k)$ acts as an orthogonal transformation
on ${\mathfrak n}_\Phi$ for each $k \in K_\Phi$,
$N_{K_\Phi}({\mathfrak n}_{\Phi,{\mathfrak v}})$
coincides with the normalizer $N_{K_\Phi}({\mathfrak v})$ of ${\mathfrak v}$
in $K_\Phi$.
Denote by $N^o_{L_\Phi}({\mathfrak n}_{\Phi,{\mathfrak v}})$ and
$N^o_{K_\Phi}({\mathfrak n}_{\Phi,{\mathfrak v}})$ the connected component
of $N_{L_\Phi}({\mathfrak n}_{\Phi,{\mathfrak v}})$ and $N_{K_\Phi}({\mathfrak n}_{\Phi,{\mathfrak v}})$
containing the identity transformation on $M$, respectively.
Then $H_{\Phi,{\mathfrak v}} = N^o_{L_\Phi}({\mathfrak n}_{\Phi,{\mathfrak v}})N_{\Phi,{\mathfrak v}}$
is a connected subgroup of $Q_\Phi$.
Assume that
$N^o_{L_\Phi}({\mathfrak n}_{\Phi,{\mathfrak v}})$ acts transitively on $F_\Phi$, which
just means that $F_\Phi \subset H_{\Phi,{\mathfrak v}} \cdot o$.
Since $H_{\Phi,{\mathfrak v}} \cap K = N^o_{L_\Phi}({\mathfrak n}_{\Phi,{\mathfrak v}}) \cap K$
and $N_{L_\Phi}({\mathfrak n}_{\Phi,{\mathfrak v}}) \cap K = N_{K_\Phi}({\mathfrak n}_{\Phi,{\mathfrak v}})
= N_{K_\Phi}({\mathfrak v})$, we see that the cohomogeneity of the action of
$H_{\Phi,{\mathfrak v}}$ on $M$ is equal to the cohomogeneity of the action
of $N^o_{K_\Phi}({\mathfrak v})$ on ${\mathfrak v}$.
Since $L_\Phi$ is reductive, we also have $N^o_{L_\Phi}({\mathfrak n}_{\Phi,{\mathfrak v}}) =
\theta N^o_{L_\Phi}({\mathfrak v})$.
Thus we get the following
construction method for cohomogeneity one actions on $M$.

\begin{prop}
\label{viagradation}
Assume that $\dim {\mathfrak n}_\Phi^1 \geq 2$ and
let $\cv$ be a subspace of
${\mathfrak n}_\Phi^1$ with $\dim {\mathfrak v} \geq 2$ such that
\begin{itemize}
\item[(i)] $N^o_{L_\Phi}({\mathfrak n}_{\Phi,{\mathfrak v}})= \theta N^o_{L_\Phi}({\mathfrak v})$
acts transitively on $F_\Phi$, and
\item[(ii)] $N^o_{K_\Phi}({\mathfrak n}_{\Phi,{\mathfrak v}}) = N^o_{K_\Phi}({\mathfrak v})$
acts transitively on the unit sphere in ${\mathfrak v}$.
\end{itemize}
Then $H_{\Phi,{\mathfrak v}} = N^o_{L_\Phi}({\mathfrak n}_{\Phi,{\mathfrak v}})N_{\Phi,{\mathfrak v}}$
acts on $M$ with cohomogeneity one
and $H_{\Phi,{\mathfrak v}} \cdot o$ is a singular orbit of this action containing $F_\Phi$.
Moreover, if ${\mathfrak v}_1$ and ${\mathfrak v}_2$ are two such subspaces which are
conjugate by an element in $K_\Phi$, then the cohomogeneity one actions of
$H_{\Phi,{\mathfrak v}_1}$ and $H_{\Phi,{\mathfrak v}_2}$ on $M$ are orbit equivalent.
\end{prop}

\begin{proof}
We only have to prove the statement about orbit equivalence. Assume that
${\rm Ad}(k)({\mathfrak v}_1) = {\mathfrak v}_2$ for some $k \in K_\Phi$. Since ${\rm Ad}(k)$
preserves the Chevalley decomposition ${\mathfrak q}_\Phi = {\mathfrak l}_\Phi \oplus {\mathfrak n}_\Phi$,
we have
${\rm Ad}(k)(N_{{\mathfrak l}_\Phi}({\mathfrak n}_{\Phi,{\mathfrak v}_1}) \oplus
{\mathfrak n}_{\Phi,{\mathfrak v}_1})
= {\rm Ad}(k)N_{{\mathfrak l}_\Phi}({\mathfrak n}_{\Phi,{\mathfrak v}_1}) \oplus
{\rm Ad}(k){\mathfrak n}_{\Phi,{\mathfrak v}_1}$.
By assumption, we have ${\rm Ad}(k)({\mathfrak v}_1) = {\mathfrak v}_2$, and since
${\rm Ad}(k)$ acts as an orthogonal transformation on ${\mathfrak n}_\Phi$, this implies
${\rm Ad}(k){\mathfrak n}_{\Phi,{\mathfrak v}_1} = {\mathfrak n}_{\Phi,{\mathfrak v}_2}$.
If $X \in N_{{\mathfrak l}_\Phi}({\mathfrak n}_{\Phi,{\mathfrak v}_1})$, we get
$$
[{\rm Ad}(k)X , {\mathfrak n}_{\Phi,{\mathfrak v}_2} ] =
[{\rm Ad}(k)X , {\rm Ad}(k){\mathfrak n}_{\Phi,{\mathfrak v}_1} ] =
{\rm Ad}(k) [ X , {\mathfrak n}_{\Phi,{\mathfrak v}_1} ] \subset
{\rm Ad}(k) {\mathfrak n}_{\Phi,{\mathfrak v}_1} = {\mathfrak n}_{\Phi,{\mathfrak v}_2},
$$
which implies
${\rm Ad}(k)N_{{\mathfrak l}_\Phi}({\mathfrak n}_{\Phi,{\mathfrak v}_1}) \subset
N_{{\mathfrak l}_\Phi}({\mathfrak n}_{\Phi,{\mathfrak v}_2})$. By an analogous argumentation we obtain
${\rm Ad}(k^{-1})N_{{\mathfrak l}_\Phi}({\mathfrak n}_{\Phi,{\mathfrak v}_2}) \subset
N_{{\mathfrak l}_\Phi}({\mathfrak n}_{\Phi,{\mathfrak v}_1})$. Altogether this shows that
${\rm Ad}(k)N_{{\mathfrak l}_\Phi}({\mathfrak n}_{\Phi,{\mathfrak v}_1}) =
N_{{\mathfrak l}_\Phi}({\mathfrak n}_{\Phi,{\mathfrak v}_2})$,
and therefore
${\rm Ad}(k){\mathfrak h}_{\Phi,{\mathfrak v}_1} = {\mathfrak h}_{\Phi,{\mathfrak v}_2}$.
Since both $H_{\Phi,{\mathfrak v}_1}$ and $H_{\Phi,{\mathfrak v}_2}$ are connected,
this implies that the actions
of $H_{\Phi,{\mathfrak v}_1}$ and $H_{\Phi,{\mathfrak v}_2}$ are orbit equivalent.
\end{proof}

We now discuss the second construction method in more detail for
maximal proper parabolic subgroups. Any such subgroup is conjugate to $Q_{\Phi_j}$ for some $j \in \{1,\ldots,r\}$,
where $\Phi_j = \Lambda \setminus \{\alpha_j\}$.
In the following we will replace the ``index $\Phi_j$'' by the ``index $j$'', that is,
the parabolic subalgebra ${\mathfrak q}_{\Phi_j}$ will be denoted by
${\mathfrak q}_j$, and so on. We discuss now a few examples of cohomogeneity one actions arising
from this construction method.

\medskip
\noindent {\it Examples.}
1. Assume that the rank of $M$ is equal to one. Thus $M$ is isometric to a hyperbolic space ${\mathbb F}H^n$
over a normed real division algebra ${\mathbb F} \in \{{\mathbb R},{\mathbb C},{\mathbb H},{\mathbb O}\}$.
In this case there is just one simple root $\alpha = \alpha_1$, and therefore $\Phi_1 = \emptyset$. The maximal
proper parabolic subgroup $Q_1$ is therefore a minimal parabolic subgroup. The parabolic subalgebra
${\mathfrak q}_1$ is given by ${\mathfrak q}_1 = {\mathfrak k}_0 \oplus {\mathfrak a} \oplus {\mathfrak n}$
with ${\mathfrak n} = {\mathfrak g}_\alpha \oplus {\mathfrak g}_{2\alpha}$. More explicitly, we have the following table:

\medskip
\begin{center}
\begin{tabular}{|l|l|l|l|l|l|}
\hline
$M$ & $G$ & $K$ & $K_0$ & ${\mathfrak g}_\alpha$ & ${\mathfrak n}$\\
\hline \hline
${\mathbb R}H^n$ & $SO^o_{1,n}$ & $SO_n$ & $SO_{n-1}$ & ${\mathbb R}^{n-1}$ & ${\mathbb R}^{n-1}$\\
${\mathbb C}H^n$ & $SU_{1,n}$ & $U_n$ & $U_{n-1}$ & ${\mathbb C}^{n-1}$ & ${\mathbb C}^{n-1} \oplus {\mathbb R}$\\
${\mathbb H}H^n$ & $Sp_{1,n}$ & $Sp_1Sp_n$ & $Sp_1Sp_{n-1}$ & ${\mathbb H}^{n-1}$ & ${\mathbb H}^{n-1} \oplus {\mathbb R}^3$\\
${\mathbb O}H^2$ & $F_4^{-20}$ & $Spin_9$ & $Spin_7$ & ${\mathbb O}$ & ${\mathbb O} \oplus {\mathbb R}^7$\\
\hline
\end{tabular}
\end{center}
\medskip

Condition (i) in Proposition \ref{viagradation} is automatically satisfied since the boundary component
$B_1$ consists of the single point $o$. Condition (ii) is equivalent to the problem:
Find all $k$-dimensional $(k \geq 2)$
linear subspaces ${\mathfrak v}$ of ${\mathfrak g}_\alpha$ for which
there exists a subgroup of $K_0$ acting transitively on the unit
sphere in ${\mathfrak v}$. The authors solved this problem in \cite{BT3} for
${\mathbb F} \in \{{\mathbb R},{\mathbb C},{\mathbb O}\}$, whereas for ${\mathbb F} = {\mathbb H}$
we only found some examples but achieved no complete classification.

If ${\mathbb F} = {\mathbb R}$, we can choose any linear subspace ${\mathfrak v} \subset {\mathbb R}^{n-1}$. However, in this case the orbit $H_{1,{\mathfrak v}} \cdot o$ is always totally geodesic in ${\mathbb R}H^n$.

If ${\mathbb F} = {\mathbb C}$, we can choose any linear subspace ${\mathfrak v} \subset {\mathbb C}^{n-1}$ with constant K\"{a}hler angle $\varphi \in [0,\pi/2]$. If $0 < \varphi < \pi/2$, then the cohomogeneity one action on ${\mathbb C}H^n$ by $H_{1,{\mathfrak v}}$ has a non-totally geodesic singular orbit and is not orbit equivalent to a cohomogeneity one action obtained by any of the other construction methods.

If ${\mathbb F} = {\mathbb O}$, we can choose any linear subspace ${\mathfrak v} \subset {\mathbb O}$ of dimension $k \in \{2,3,4,6,7\}$. The cohomogeneity one action
on ${\mathbb O}H^2$ by $H_{1,{\mathfrak v}}$ has a non-totally geodesic singular orbit and is not
orbit equivalent to a cohomogeneity one action obtained by any of the other construction methods.

If ${\mathbb F} ={\mathbb H}$, we can choose linear subspaces ${\mathfrak v} \subset {\mathbb H}^{n-1}$ with constant quaternionic K\"{a}hler angle. However, the classification of such subspaces is not yet finalized.

\smallskip
2. Let $M = G_2^2/SO_4$. Then $\dim M = 8 $, and ${\mathfrak g} =
{\mathfrak g}_2^2$ is a split real form of ${\mathfrak g}^{\mathbb
C}$.
For $M = G_2^2/SO_4$ the corresponding root system $\Sigma$ is of type $(G_2)$ and
all root spaces have real dimension $1$. We label the simple roots
by $\alpha_1$ and $\alpha_2$ so that $3\alpha_1 + 2\alpha_2$ is the
highest root in $\Sigma^+$, and choose $j = 1$, that is, $\Phi_1 = \{\alpha_2\}$. Then we have
${\mathfrak k}_0 = \{0\}$ and ${\mathfrak g}_0 = {\mathfrak a} \cong {\mathbb R}^2$. Moreover,
\begin{eqnarray*}
{\mathfrak n}_1^1 & = & {\mathfrak g}_{\alpha_1} \oplus {\mathfrak g}_{\alpha_1 +
\alpha_2} \cong {\mathbb R}^2,\\
{\mathfrak n}_1^2 & = & {\mathfrak g}_{2\alpha_1 + \alpha_2} \cong {\mathbb R},\\
{\mathfrak n}_1^3 & = & {\mathfrak g}_{3\alpha_1+\alpha_2} \oplus {\mathfrak g}_{3\alpha_1
+ 2\alpha_2} \cong {\mathbb R}^2,\\
{\mathfrak n}_1 = {\mathfrak n}_1^1 \oplus {\mathfrak n}_1^2 \oplus {\mathfrak n}_1^3 & = & {\mathfrak n} =
{\mathfrak g}_{\alpha_1} \oplus {\mathfrak g}_{\alpha_1 +
\alpha_2} \oplus {\mathfrak g}_{2\alpha_1 + \alpha_2} \oplus
{\mathfrak g}_{3\alpha_1+\alpha_2} \oplus {\mathfrak g}_{3\alpha_1
+ 2\alpha_2} \cong {\mathbb R}^5,\\
{\mathfrak a}_1 & = & {\rm ker}\,\alpha_2 = {\mathbb R}H^1 \cong {\mathbb R},\\
{\mathfrak m}_1 = {\mathfrak g}_1 & = & {\mathfrak g}_{-\alpha_2} \oplus {\mathbb R}H^2 \oplus
{\mathfrak g}_{\alpha_2} \cong {\mathfrak s}{\mathfrak l}_2({\mathbb R}),\\
{\mathfrak l}_1 = {\mathfrak m}_1 \oplus {\mathfrak a}_1
& = & {\mathfrak g}_{-\alpha_2} \oplus {\mathfrak g}_0 \oplus
{\mathfrak g}_{\alpha_2} \cong {\mathfrak s}{\mathfrak l}_2({\mathbb
R}) \oplus {\mathbb R},\\
{\mathfrak k}_1 & = & {\mathfrak k}_{\alpha_2} \cong {\mathfrak s}{\mathfrak o}_2.
\end{eqnarray*}
This explicit description shows that $F_1 \cong SL_2({\mathbb R})/SO_2 \times {\mathbb E} =
{\mathbb R}H^2 \times {\mathbb E}$ and
that $K_1^o \cong SO_2$ acts transitively on the unit sphere in
${\mathfrak v} = {\mathfrak n}_1^1 \cong {\mathbb R}^2$.
It follows that $H_{1,{\mathfrak v}}$ acts on $M$ with cohomogeneity one
whose singular orbit has codimension $2$ and contains $F_1 \cong
{\mathbb R}H^2 \times {\mathbb E}$. The Lie algebra of $H_{1,{\mathfrak v}}$
is given by
$$
{\mathfrak h}_{1,{\mathfrak v}} = {\mathfrak g}_{-\alpha_2} \oplus
{\mathfrak g}_0 \oplus {\mathfrak g}_{\alpha_2}
\oplus {\mathfrak g}_{2\alpha_1 + \alpha_2} \oplus
{\mathfrak g}_{3\alpha_1+\alpha_2} \oplus {\mathfrak g}_{3\alpha_1
+ 2\alpha_2}.
$$

\smallskip
3. Let $M = G_2^{\mathbb C}/G_2$. Then $\dim M = 14$, and
the corresponding root system $\Sigma$ is of type $(G_2)$ and can be identified with the
root system of the complex simple Lie algebra $({\mathfrak g}_2)^{\mathbb C}$. Therefore all root spaces
have complex dimension $1$. As in the previous example we label the simple roots
by $\alpha_1$ and $\alpha_2$ so that $3\alpha_1 + 2\alpha_2$ is the
highest root in $\Sigma^+$, and choose again $j = 1$, and hence $\Phi_1 = \{\alpha_2\}$. Then we have
${\mathfrak k}_0 \cong {\mathfrak u}_1 \oplus {\mathfrak u}_1$,
${\mathfrak g}_0 \cong {\mathfrak u}_1 \oplus {\mathfrak u}_1 \oplus {\mathbb R}^2$, and
\begin{eqnarray*}
{\mathfrak n}_1^1 & = & {\mathfrak g}_{\alpha_1} \oplus {\mathfrak g}_{\alpha_1 +
\alpha_2} \cong {\mathbb C}^2,\\
{\mathfrak n}_1^2 & = & {\mathfrak g}_{2\alpha_1 + \alpha_2} \cong {\mathbb C},\\
{\mathfrak n}_1^3 & = & {\mathfrak g}_{3\alpha_1+\alpha_2} \oplus {\mathfrak g}_{3\alpha_1
+ 2\alpha_2} \cong {\mathbb C}^2,\\
{\mathfrak n}_1 = {\mathfrak n}_1^1 \oplus {\mathfrak n}_1^2 \oplus {\mathfrak n}_1^3 & = & {\mathfrak n} =
{\mathfrak g}_{\alpha_1} \oplus {\mathfrak g}_{\alpha_1 +
\alpha_2} \oplus {\mathfrak g}_{2\alpha_1 + \alpha_2} \oplus
{\mathfrak g}_{3\alpha_1+\alpha_2} \oplus {\mathfrak g}_{3\alpha_1
+ 2\alpha_2} \cong {\mathbb C}^5,\\
{\mathfrak l}_1 & = & {\mathfrak g}_{-\alpha_2} \oplus {\mathfrak g}_0 \oplus
{\mathfrak g}_{\alpha_2} \cong {\mathfrak s}{\mathfrak l}_2({\mathbb C})
\oplus {\mathfrak u}_1 \oplus {\mathbb R},\\
{\mathfrak a}_1 & = & {\rm ker}\,\alpha_2 = {\mathbb R}H^1 \cong {\mathbb R},\\
{\mathfrak m}_1 = {\mathfrak l}_1 \ominus {\mathfrak a}_1 & \cong &
{\mathfrak s}{\mathfrak l}_2({\mathbb C}) \oplus {\mathfrak u}_1,\\
{\mathfrak g}_1 = [{\mathfrak m}_1,{\mathfrak m}_1] & \cong & {\mathfrak s}{\mathfrak l}_2({\mathbb C}),\\
{\mathfrak z}_1 & \cong & {\mathfrak u}_1,\\
{\mathfrak k}_1 & = & {\mathfrak k}_{\alpha_2} \oplus {\mathfrak u}_1 \cong {\mathfrak s}{\mathfrak u}_2
\oplus {\mathfrak u}_1 \cong {\mathfrak u}_2.
\end{eqnarray*}
From this we see that $F_1 \cong SL_2({\mathbb C})/SU_2 \times {\mathbb E} \cong {\mathbb R}H^3 \times {\mathbb E}$.
Moreover, $K_1^o \cong U_2$ acts transitively on the unit
sphere in ${\mathfrak v} = {\mathfrak n}_1^1 \cong {\mathbb C}^2$.
It follows that $H_{1,{\mathfrak v}}$ acts on $M$ with cohomogeneity one
whose singular orbit has codimension $4$ and contains $F_1 \cong
{\mathbb R}H^3 \times {\mathbb E}$. The Lie algebra of $H_{1,{\mathfrak v}}$
is given by
$$
{\mathfrak h}_{1,{\mathfrak v}} = {\mathfrak g}_{-\alpha_2} \oplus {\mathfrak g}_0 \oplus {\mathfrak g}_{\alpha_2}
\oplus {\mathfrak g}_{2\alpha_1 + \alpha_2} \oplus
{\mathfrak g}_{3\alpha_1+\alpha_2} \oplus {\mathfrak g}_{3\alpha_1
+ 2\alpha_2}.
$$

\smallskip
4. The following example illustrates that the two different construction methods can lead to
orbit equivalent cohomogeneity one actions even when $|\Lambda \setminus \Phi| = 1$.
Let $M = SO^o_{2,n+2}/SO_2SO_{n+2}$ and $n \geq 1$. Then $\dim M = 2n+4$ and
the corresponding root system $\Sigma$ is of type $(B_2)$. Let $\alpha_1$ and $\alpha_2$
be corresponding simple roots such that $\alpha_1$ is the longer of the two roots. Then
we have $\Sigma^+ = \{\alpha_1,\alpha_2,\alpha_1+\alpha_2,\alpha_1+2\alpha_2\}$, and the multiplicities
of the two long roots $\alpha_1$ and $\alpha_1 + 2\alpha_2$ are $1$
and of the two short roots $\alpha_2$ and $\alpha_1 + \alpha_2$ are $n$.
We have ${\mathfrak k}_0 \cong {\mathfrak s}{\mathfrak o}_n$ and ${\mathfrak a} \cong {\mathbb R}^2$.

Firstly, we choose $j = 1$, that is, $\Phi_1 = \{\alpha_2\}$. Then we have
\begin{eqnarray*}
{\mathfrak n}_1 = {\mathfrak n}_1^1 & = & {\mathfrak g}_{\alpha_1} \oplus {\mathfrak g}_{\alpha_1 +
\alpha_2} \oplus {\mathfrak g}_{\alpha_1 + 2\alpha_2} \cong {\mathbb R}^{n+2},\\
{\mathfrak l}_1 & = & {\mathfrak g}_{-\alpha_2} \oplus {\mathfrak g}_0 \oplus
{\mathfrak g}_{\alpha_2} \cong {\mathfrak s}{\mathfrak o}_{1,n+1} \oplus {\mathbb R},\\
{\mathfrak a}_1 & = & {\rm ker}\,\alpha_2 = {\mathbb R}H^1 \cong {\mathbb R},\\
{\mathfrak g}_1 = {\mathfrak m}_1 = {\mathfrak l}_1 \ominus {\mathfrak a}_1 & \cong &
{\mathfrak s}{\mathfrak o}_{1,n+1}\\
{\mathfrak k}_1 & = & {\mathfrak s}{\mathfrak o}_{n+1}.
\end{eqnarray*}
From this we see that $F_1 \cong SO^o_{1,n+1}/SO_{n+1} \times {\mathbb E}$ and
$K_1^o \cong SO_{n+1}$. The $K_1^o$-module ${\mathfrak n}_1$ decomposes into a
$1$-dimensional trivial module
${\mathfrak n}_{1,0} \cong {\mathbb R} \subset {\mathfrak g}_{\alpha_1} \oplus
{\mathfrak g}_{\alpha_1 + 2\alpha_2} \cong {\mathbb R}^2$
and an irreducible module
${\mathfrak v} \cong {\mathbb R}^{n+1} \supset {\mathfrak g}_{\alpha_1 + \alpha_2} \cong {\mathbb R}^n$.
The action of $K_1^o \cong SO_{n+1}$ on the irreducible module
${\mathfrak v} \cong {\mathbb R}^{n+1}$ is the standard one and acts
transitively on the unit sphere.
It follows that $H_{1,{\mathfrak v}}$ acts on $M$ with cohomogeneity one
whose singular orbit $W$ has codimension $n+1$ and contains $F_1 \cong
{\mathbb R}H^{n+1} \times {\mathbb E}$. The Lie algebra of $H_{1,{\mathfrak v}}$
is given by
$$
{\mathfrak h}_{1,{\mathfrak v}} = {\mathfrak g}_{-\alpha_2} \oplus {\mathfrak g}_0 \oplus {\mathfrak g}_{\alpha_2}
\oplus {\mathfrak n}_{1,0} = {\mathfrak g}_1 \oplus ({\mathfrak a}_1 \oplus {\mathfrak n}_{1,0}).
$$
However, the orbit through $o$ of the connected subgroup of $G$
with Lie algebra ${\mathfrak a}_1 \oplus {\mathfrak n}_{1,0}$ is a totally geodesic real hyperbolic
plane ${\mathbb R}H^2$, and hence $W$ is a totally geodesic submanifold of $M$
which is congruent to the Riemannian product ${\mathbb R}H^{n+1} \times {\mathbb R}H^2$
of two real hyperbolic spaces.

Finally, we choose $j = 2$, that is, $\Phi_2 = \{\alpha_1\}$. Then we have
\begin{eqnarray*}
{\mathfrak n}_2^1 & = & {\mathfrak g}_{\alpha_2} \oplus {\mathfrak g}_{\alpha_1 +
\alpha_2} \cong {\mathbb R}^{2n},\\
{\mathfrak n}_2^2 & = & {\mathfrak g}_{\alpha_1 + 2\alpha_2} \cong {\mathbb R}\\
{\mathfrak n}_2 & = & {\mathfrak g}_{\alpha_2} \oplus {\mathfrak g}_{\alpha_1 +
\alpha_2} \oplus {\mathfrak g}_{\alpha_1 + 2\alpha_2} \cong {\mathbb R}^{2n+1}\\
{\mathfrak l}_2 & = & {\mathfrak g}_{-\alpha_1} \oplus {\mathfrak g}_0 \oplus
{\mathfrak g}_{\alpha_1} \cong {\mathfrak s}{\mathfrak o}_{1,2} \oplus {\mathfrak s}{\mathfrak o}_n \oplus {\mathbb R},\\
{\mathfrak a}_2 & = & {\rm ker}\,\alpha_1 = {\mathbb R}H^2 \cong {\mathbb R},\\
{\mathfrak g}_2 = {\mathfrak m}_2 = {\mathfrak l}_2 \ominus {\mathfrak a}_2 & \cong &
 {\mathfrak s}{\mathfrak o}_{1,2} \oplus {\mathfrak s}{\mathfrak o}_n\\
{\mathfrak k}_2 & \cong &  {\mathfrak s}{\mathfrak o}_2 \oplus {\mathfrak s}{\mathfrak o}_n.
\end{eqnarray*}
From this we see that $F_2 \cong SO^o_{1,2}/SO_2 \times {\mathbb E} = {\mathbb R}H^2 \times {\mathbb E}$ and $K_2^o \cong SO_2SO_n$.
The representation of $K_2^o$ on ${\mathfrak n}_2^1$ is isomorphic to the tensor representation
of $SO_2SO_n$ on ${\mathbb R}^2 \otimes {\mathbb R}^n \cong {\mathbb R}^{2n}$.

The symmetric space $M = SO^o_{2,n+2}/SO_2SO_{n+2}$ is Hermitian and hence has a natural complex structure $J$.
This complex structure turns ${\mathfrak n}_2^1 \cong {\mathbb R}^{2n}$ into a complex
vector space ${\mathbb C}^n$ so that ${\mathfrak g}_{\alpha_2}$ and ${\mathfrak g}_{\alpha_1 +
\alpha_2}$ are real subspaces which are mapped onto each
other by $J$. Moreover, the action of $SO_2 \subset SO_2SO_n \cong K_2^o$ on ${\mathfrak n}_2^1$ is isomorphic to
the standard action of the circle group on ${\mathbb C}^n$, and the action of
$SO_n \subset SO_2SO_n \cong K_2^o$ on ${\mathbb R}^n \cong {\mathfrak g}_{\alpha_2} \subset {\mathfrak n}_2^1$
and on ${\mathbb R}^n \cong {\mathfrak g}_{\alpha_1 + \alpha_2} \subset {\mathfrak n}_2^1$
is isomorphic to the standard action of $SO_n$ on ${\mathbb R}^n$. We now construct cohomogeneity one
actions on $M$ through two different types of subspaces ${\mathfrak v}$.

Firstly, let ${\mathfrak v}$ be a $k$-dimensional linear subspace of ${\mathfrak g}_{\alpha_2}$ with
$k \geq 2$. Then $N^o_{K_2}({\mathfrak v})$ is isomorphic to
$SO_kSO_{n-k} \subset SO_n \subset SO_2SO_n$ and acts transitively on the unit sphere
in ${\mathfrak v}$. Moreover,
$$N_{{\mathfrak l}_2}({\mathfrak n}_{2,{\mathfrak v}}) =
({\mathfrak s}{\mathfrak o}_k \oplus {\mathfrak s}{\mathfrak o}_{n-k}) \oplus
{\mathfrak a} \oplus {\mathfrak g}_{\alpha_1}$$
where ${\mathfrak s}{\mathfrak o}_k \oplus {\mathfrak s}{\mathfrak o}_{n-k} \cong N_{{\mathfrak k}_2}({\mathfrak v})$.
We easily see that the connected subgroup of $G$ with Lie algebra ${\mathfrak a} \oplus {\mathfrak g}_{\alpha_1}$
acts transitively on $F_2 \cong SO^o_{1,2}/SO_2 \times {\mathbb E}$.
Altogether it follows that $H_{2,{\mathfrak v}}$ acts on $M$ with cohomogeneity one and
with a singular orbit $W$ of codimension $k$ and containing $F_2$. The Lie algebra of $H_{2,{\mathfrak v}}$
is given by
$$
{\mathfrak h}_{2,{\mathfrak v}} =  ({\mathfrak s}{\mathfrak o}_k \oplus {\mathfrak s}{\mathfrak o}_{n-k})
 \oplus {\mathfrak a} \oplus {\mathfrak g}_{\alpha_1}
\oplus ({\mathfrak g}_{\alpha_2} \ominus {\mathfrak v})\oplus {\mathfrak g}_{\alpha_1 +
\alpha_2} \oplus {\mathfrak g}_{\alpha_1 + 2\alpha_2}.
$$
Note that ${\mathfrak s}{\mathfrak o}_k \oplus {\mathfrak
s}{\mathfrak o}_{n-k} \subset {\mathfrak s}{\mathfrak o}_n \cong
{\mathfrak k}_0$. However, it is evident from the explicit description of
${\mathfrak h}_{2,{\mathfrak v}}$ that the action of $H_{2,{\mathfrak v}}$ on $M$ is orbit
equivalent to the action of the canonical extension of a cohomogeneity one action on
the boundary component $B_1 = SO^o_{1,n+1}/SO_{n+1}$. Instead of picking a subspace
${\mathfrak v}$ of the real subspace ${\mathfrak g}_{\alpha_2}$ of
${\mathfrak n}_2^1$, we could also select a subspace ${\mathfrak v}$
of any of the real subspaces of ${\mathfrak n}_2^1$ obtained by
rotating ${\mathfrak g}_{\alpha_2}$ in ${\mathfrak n}_2^1$ by means
of the $SO_2$-action with the $SO_2$ whose Lie algebra is
${\mathfrak s}{\mathfrak o}_2$ in ${\mathfrak k}_2 \cong {\mathfrak
s}{\mathfrak o}_2 \oplus {\mathfrak s}{\mathfrak o}_n$. For example,
${\mathfrak g}_{\alpha_1+\alpha_2}$ is such a subspace. However,
such a cohomogeneity one action is conjugate to one constructed from
a subspace in ${\mathfrak g}_{\alpha_2}$.

\section{Proof of Theorem \ref{maintheorem}}
\label{section:classification}

Let $H$ be a connected subgroup of $G$ acting on
$M$ with cohomogeneity one. If the orbits form a Riemannian foliation, a complete classification
up to orbit equivalence was obtained by the authors in \cite{BT1} for irreducible
symmetric spaces $M$. For reducible symmetric spaces the corresponding problem is still unsolved.
We assume from now on that the action has a singular orbit $W$.
Then $H$ is contained either in a proper maximal reductive subgroup of $G$
or in a proper maximal parabolic subgroup of $G$.
In the first case we have a totally geodesic singular orbit (see Theorem \ref{maximal}).
For irreducible symmetric spaces $M$ the classification of such actions was obtained by the authors
in \cite{BT2}. For reducible symmetric spaces the corresponding problem is still not solved.
We assume from now on that $H$
is contained in a proper maximal parabolic subgroup of $G$, or equivalently,
${\mathfrak h} \subset {\mathfrak q}_j$ for some $j \in \{1,\ldots,r\}$. Without loss of generality we
may assume that $o \in W$, that is $W = H \cdot o$.

Consider the slice
representation
$$
\chi : H\cap K \to O(\nu_oW),\ k \mapsto d_ok|_{\nu_oW}.
$$
Since $H \subset Q_j$, we have
$ H \cap K \subset Q_j \cap K = K_j$, and therefore
$d_ok(\xi) = {\rm Ad}(k)\xi$ for all $\xi \in \nu_oW$ and $k \in H\cap K$, where we
identify
$$
T_oM \cong ({\mathfrak l}_j \cap {\mathfrak p}) \oplus {\mathfrak n}_j
= {\mathfrak b}_j \oplus {\mathfrak a}_j
\oplus \left( \bigoplus_{\nu = 1}^{m_j} {\mathfrak n}_j^\nu \right).
$$
Recall that $T_oB_j \cong {\mathfrak b}_j$ under the above identification.

We first show that
the normal space $\nu_o W$ is contained in either
${\mathfrak b}_j \cong T_oB_j$ or ${\mathfrak n}_j^1$.
First of all, we use the fact that
$H \cap K = H \cap K_j$ acts transitively on the unit sphere in $\nu_o W$.
We decompose the parabolic subalgebra ${\mathfrak q}_j$ into
 ${\mathfrak q}_j = {\mathfrak k}_j \oplus {\mathfrak b}_j \oplus
{\mathfrak a}_j \oplus {\mathfrak n}_j$ and denote by
$\tau : {\mathfrak q}_j \to {\mathfrak b}_j \oplus {\mathfrak a}_j \oplus  {\mathfrak n}_j$
the canonical projection with respect to this decomposition. Then we have
$\tau({\mathfrak h}) = T_oW = ({\mathfrak b}_j \oplus {\mathfrak a}_j \oplus  {\mathfrak n}_j)
\ominus \nu_oW$.
Since $\dim {\mathfrak a}_j = 1$, we must have
\begin{eqnarray}
\label{eq:a-tangent}
\nu_o W \subset {\mathfrak b}_j \oplus {\mathfrak n}_j .
\end{eqnarray}
Let us define
\begin{eqnarray*}
(\nu_o W)_0
& := &
{\mathfrak b}_j \ominus (\tau({\mathfrak h}) \cap {\mathfrak b}_j) , \\
(\nu_o W)_\nu
& := &
{\mathfrak n}_j^\nu \ominus (\tau({\mathfrak h}) \cap {\mathfrak n}_j^\nu)
\quad \mbox{for $\nu = 1 , \ldots , m_j$.}
\end{eqnarray*}
It is easy to see that
\begin{eqnarray}
\label{eq:nu_0-dec}
\nu_o W \subset (\nu_o W)_0 \oplus (\nu_o W)_1 \oplus \cdots \oplus (\nu_o W)_{m_j} .
\end{eqnarray}

\begin{lm}
\label{lm:irr-module}
Let $0 \neq X \in \nu_o W$, and denote by
$\pi_k : (\nu_o W)_0 \oplus \cdots \oplus (\nu_o W)_{m_j} \to (\nu_o W)_k$
the canonical projection.
Then we have
\begin{itemize}
\item[(1)]
$\nu_o W = \R X \oplus [ {\mathfrak h} \cap {\mathfrak k}_j , X ]$.
\item[(2)]
$\pi_k |_{\nu_o W} : \nu_o W \to (\nu_o W)_k$ is onto and
$(H \cap K_j)$-equivariant.
\item[(3)]
If $\pi_k(X) = 0$ then $(\nu_o W)_k = 0$.
\item[(4)]
If $(\nu_o W)_k \neq 0$ then
$\pi_k |_{\nu_o W} : \nu_o W \to (\nu_o W)_k$ is an $(H \cap K_j)$-equivariant isomorphism.
In particular, $(\nu_o W)_k$ is an irreducible $({\mathfrak h} \cap {\mathfrak k}_j)$-module.
\end{itemize}
\end{lm}

\begin{proof}
Since $H \cap K_j$ acts transitively on the unit sphere in $\nu_o W$,
the subspace
$[ {\mathfrak h} \cap {\mathfrak k}_j , X ]$
in $\nu_o W$ has codimension one.
This implies (1) since
$[ {\mathfrak h} \cap {\mathfrak k}_j , X ]$
is perpendicular to $X$.
Statement (2) follows from the fact that
$H \cap K_j$ preserves the decomposition (\ref{eq:nu_0-dec}).
To show (3), assume that $\pi_k(X) = 0$.
This means $\langle X , (\nu_o W)_k \rangle = 0$.
Since ${\mathfrak h} \cap {\mathfrak k}_j$ normalizes $(\nu_o W)_k$ and by (1),
we have $\langle \nu_o W , (\nu_o W)_k \rangle = 0$.
This implies $(\nu_o W)_k \subset \tau({\mathfrak h})$ and
hence $(\nu_o W)_k = 0$.
To show (4), assume that $(\nu_o W)_k \neq 0$.
Then $\pi_k |_{\nu_o W}$ is injective by (3), and taking into account (2),
we see that $\pi_k |_{\nu_o W}$ is an isomorphism.
\end{proof}

In the second step
we use ${\mathfrak a}_j \subset \tau({\mathfrak h})$,
which follows directly from (\ref{eq:a-tangent}).
By definition,
there exists $H^j_{{\mathfrak k}} \in {\mathfrak k}_j$ such that
$H^j_{{\mathfrak k}} + H^j \in {\mathfrak h}$.
Note that
$H^j_{{\mathfrak k}} \in ({\mathfrak h})_{{\mathfrak k}_j}$,
where $({\mathfrak h})_{{\mathfrak k}_j}$ is obtained by orthogonally
projecting ${\mathfrak h}$ into ${\mathfrak k}_j$.
We decompose this subspace orthogonally into
\begin{eqnarray}
\label{eq:(h)_k-dec}
({\mathfrak h})_{{\mathfrak k}_j} = ({\mathfrak h} \cap {\mathfrak k}_j) \oplus
(({\mathfrak h})_{{\mathfrak k}_j} \ominus ({\mathfrak h} \cap {\mathfrak k}_j)) .
\end{eqnarray}
If we write
$H^j_{{\mathfrak k}} = (H^j_{{\mathfrak k}})_1 + (H^j_{{\mathfrak k}})_2$
according to this decomposition,
then
$(H^j_{{\mathfrak k}})_1 \in {\mathfrak h} \cap {\mathfrak k}_j \subset {\mathfrak h}$
and hence
$(H^j_{{\mathfrak k}})_2 + H^j \in {\mathfrak h}$.
By this argument we may and do assume that
$$
H^j_{{\mathfrak k}} \in
({\mathfrak h})_{{\mathfrak k}_j} \ominus ({\mathfrak h} \cap {\mathfrak k}_j) .
$$
In the next lemma we investigate the action of
$$
f := \ad(H^j_{{\mathfrak k}}) .
$$

\begin{lm}
\label{lm:f^2}
For each $k = 0, 1, \ldots, m_j$, we have
\begin{itemize}
\item[(1)]
$f$ normalizes $(\nu_o W)_k$,
\item[(2)]
$f^2 = -c_k^2 \cdot \id$ on $(\nu_o W)_k$.
\end{itemize}
\end{lm}

\begin{proof}
We first show (1).
Since $H^j_{{\mathfrak k}} \in {\mathfrak k}$,
the map $f$ is skewsymmetric.
Therefore it is enough to show that $f$ normalizes
$\tau({\mathfrak h}) \cap {\mathfrak b}_j$
and $\tau({\mathfrak h}) \cap {\mathfrak n}_j^\nu$.
Let $X \in \tau({\mathfrak h}) \cap {\mathfrak b}_j$.
There exists $X_{{\mathfrak k}} \in {\mathfrak k}_j$ such that
$X_{{\mathfrak k}} + X \in {\mathfrak h}$.
Since $[ {\mathfrak a}_j , {\mathfrak m}_j] = 0$,
we have
$$
{\mathfrak h}
\ni [ H^j_{{\mathfrak k}} + H^j , X_{{\mathfrak k}} + X ]
= [ H^j_{{\mathfrak k}} , X_{{\mathfrak k}} ]
+ [ H^j_{{\mathfrak k}} , X ] .
$$
This concludes
$f(X) = [ H^j_{{\mathfrak k}} , X ] \in \tau({\mathfrak h}) \cap {\mathfrak b}_j$.
Next, let $Y \in \tau({\mathfrak h}) \cap {\mathfrak n}_j^\nu$.
There exists $Y_{{\mathfrak k}} \in {\mathfrak k}_j$ such that
$Y_{{\mathfrak k}} + Y \in {\mathfrak h}$.
By definition of $H^j$, we have
$$
{\mathfrak h}
\ni [ H^j_{{\mathfrak k}} + H^j , Y_{{\mathfrak k}} + Y ]
= [ H^j_{{\mathfrak k}} , Y_{{\mathfrak k}} ]
+ [ H^j_{{\mathfrak k}} , Y ] + \nu Y .
$$
Hence we have
$[ H^j_{{\mathfrak k}} , Y ] + \nu Y \in \tau({\mathfrak h}) \cap {\mathfrak n}_j^\nu$.
Since $Y \in \tau({\mathfrak h}) \cap {\mathfrak n}_j^\nu$ by assumption,
we conclude that
$f(Y) = [ H^j_{{\mathfrak k}} , Y ]
\in \tau({\mathfrak h}) \cap {\mathfrak n}_j^\nu$.
This finishes (1).

To show (2), we need
\begin{eqnarray}
\label{eq:h^j-commute}
[ {\mathfrak h} \cap {\mathfrak k}_j , H^j_{{\mathfrak k}} ] = 0 .
\end{eqnarray}
Let $X \in {\mathfrak h} \cap {\mathfrak k}_j$.
Then we have
$
[ X , H^j_{{\mathfrak k}} + H^j ]
= [ X , H^j_{{\mathfrak k}} ] \in {\mathfrak h} \cap {\mathfrak k}_j
$. On the other hand,
since ${\mathfrak h} \cap {\mathfrak k}_j$ preserves the decomposition
(\ref{eq:(h)_k-dec}),
we also have
$
[ X , H^j_{{\mathfrak k}} ] \in
({\mathfrak h})_{{\mathfrak k}_j} \ominus ({\mathfrak h} \cap {\mathfrak k}_j)
$. This implies
$[ X , H^j_{{\mathfrak k}} ] = 0$,
which finishes the proof of (\ref{eq:h^j-commute}).

We now prove (2).
By Lemma \ref{lm:irr-module} (4),
each $(\nu_o W)_k$ is an irreducible $({\mathfrak h} \cap {\mathfrak k}_j)$-module.
Hence (\ref{eq:h^j-commute}) and Schur's Lemma yield that
$f$ is a multiple of the identity on the complexification of $(\nu_o W)_k$.
Since all eigenvalues of $f$ are purely imaginary,
we conclude that $f^2 = -c_k^2 \cdot \id$ on $(\nu_o W)_k$.
\end{proof}

In the third step,
we use the fact that ${\mathfrak h}$ is a subalgebra,
and prove that $\nu_o W$ is in the suitable position.

\begin{prop}
\label{normalspace}
We have $\nu_oW \subset {\mathfrak b}_j \cong T_oB_j$
or $\nu_oW \subset {\mathfrak n}_j^1$.
\end{prop}

\begin{proof}
First we assume that $(\nu_o W)_1 = 0$,
and show that $\nu_o W \subset {\mathfrak b}_j$.
By assumption, we have ${\mathfrak n}_j^1 \subset \tau({\mathfrak h})$.
Then, for each $X,Y \in {\mathfrak n}_j^1$,
there exist $X_{{\mathfrak k}}, Y_{{\mathfrak k}} \in {\mathfrak k}_j$
such that
$X_{{\mathfrak k}} + X , Y_{{\mathfrak k}} + Y \in {\mathfrak h}$.
Since ${\mathfrak h}$ is a subalgebra, we have
$$
\tau({\mathfrak h}) \ni \tau([ X_{{\mathfrak k}} + X , Y_{{\mathfrak k}} + Y ])
= [ X_{{\mathfrak k}} , Y ] + [ X , Y_{{\mathfrak k}} ] + [ X , Y ] .
$$
Since ${\mathfrak k}_j$ normalizes ${\mathfrak n}_j^1$
we get
$[ X_{{\mathfrak k}} , Y ] + [ X , Y_{{\mathfrak k}} ] \in {\mathfrak n}_j^1
\subset \tau({\mathfrak h})$
and therefore
$[ X , Y ] \in \tau({\mathfrak h})$.
This means ${\mathfrak n}_j^2 \subset \tau({\mathfrak h})$,
since ${\mathfrak n}_j^2$ is generated by ${\mathfrak n}_j^1$.
Recall that ${\mathfrak n}_j$ is generated by ${\mathfrak n}_j^1$.
Hence, using this argument inductively,
we conclude that
${\mathfrak n}_j \subset \tau({\mathfrak h})$.
This finishes the first case.

We next assume that $(\nu_o W)_1 \neq 0$,
which is the second case.
We show that $\nu_o W \subset {\mathfrak n}_j^1$,
that is,
\begin{eqnarray}
\label{eq:nu>1}
(\nu_o W)_\nu = 0 \quad (\mbox{for $\nu \neq 1$}) .
\end{eqnarray}
Assume that
$(\nu_o W)_\nu \neq 0$ for some $\nu \neq 1$.
Let $0 \neq X = X_0 + X_1 + \cdots + X_{m_j} \in \nu_o W$,
where $X_k \in (\nu_o W)_k$.
We have $X_1 \neq 0 \neq X_\nu$ by assumption and Lemma \ref{lm:irr-module} (3).
We put
$$
Y_1 := || X_\nu ||^2 X_1 , \quad
Y_\nu := - || X_1 ||^2 X_\nu .
$$
Since $\langle X , Y_1 + Y_\nu \rangle = 0$, Lemma \ref{lm:irr-module} (1) and the skewsymmetry of ${\rm ad}(H)$ for all $H \in {\mathfrak h} \cap {\mathfrak k}_j$ imply
\begin{eqnarray}
\label{eq:eq-nu-1}
Y_1 + Y_\nu \in \tau({\mathfrak h}) .
\end{eqnarray}
There exists $Y_{{\mathfrak k}} \in {\mathfrak k}_j$
such that $Y_{{\mathfrak k}} + Y_1 + Y_\nu \in {\mathfrak h}$.
This yields
\begin{eqnarray}
\label{eq:eq-nu-2}
\tau({\mathfrak h}) \ni
\tau([ H^j_{{\mathfrak k}} + H^j , Y_{{\mathfrak k}} + Y_1 + Y_\nu ])
= (f(Y_1) + Y_1) + (f(Y_\nu) + \nu Y_\nu) .
\end{eqnarray}
By bracketing again  we get
\begin{eqnarray*}
\tau({\mathfrak h}) & \ni &
\tau([ H^j_{{\mathfrak k}} + H^j , [ H^j_{{\mathfrak k}} + H^j , Y_{{\mathfrak k}} + Y_1 + Y_\nu ]]) \\
& = &
(f^2(Y_1) + 2 f(Y_1) + Y_1) + (f^2(Y_\nu) + 2 \nu f(Y_\nu) + \nu^2 Y_\nu) .
\end{eqnarray*}
From Lemma \ref{lm:f^2}
we know that $f^2(Y_1) = - c_1^2 Y_1$ and $f^2(Y_\nu) = - c_\nu^2 Y_\nu$,
and therefore
\begin{eqnarray}
\label{eq:eq-nu-3}
((1 - c_1^2) Y_1 + 2 f(Y_1)) + ((\nu^2 - c_\nu^2) Y_\nu + 2 \nu f(Y_\nu))
\in \tau({\mathfrak h}) .
\end{eqnarray}
From (\ref{eq:eq-nu-3}) and (\ref{eq:eq-nu-2}) we get
\begin{eqnarray*}
(-1 - c_1^2) Y_1 + (\nu^2 - c_\nu^2 - 2 \nu) Y_\nu + 2 (\nu -1) f(Y_\nu)
\in \tau({\mathfrak h}) .
\end{eqnarray*}
This and (\ref{eq:eq-nu-1}) yield
\begin{eqnarray*}
Y'_\nu :=
(\nu^2 - c_\nu^2 - 2 \nu + 1 + c_1^2) Y_\nu + 2 (\nu -1) f(Y_\nu)
\in \tau({\mathfrak h}) .
\end{eqnarray*}
Therefore,
$Y'_\nu \in \tau({\mathfrak h}) \cap {\mathfrak b}_j$ if $\nu = 0$,
and
$Y'_\nu \in \tau({\mathfrak h}) \cap {\mathfrak n}_j^\nu$ if $\nu \geq 2$.
On the other hand, we have
$Y_\nu \in (\nu_o W)_\nu$ by assumption,
and $f(Y_\nu) \in (\nu_o W)_\nu$ by Lemma \ref{lm:f^2} (1).
This means $Y'_\nu \in (\nu_o W)_\nu$, and we thus get
$Y'_\nu = 0 $.
Since $\langle Y_\nu , f(Y_\nu) \rangle = 0$ and $\nu \neq 1$,
we have $f(Y_\nu) = 0$.
But this implies $c_\nu = 0$ and hence
$\nu^2 - c_\nu^2 - 2 \nu + 1 + c_1^2 > 0$,
which contradicts $Y'_\nu = 0$.
This finishes our claim (\ref{eq:nu>1})
and shows that $\nu_o W \subset {\mathfrak n}_j^1$.
\end{proof}

We next study the structure of ${\mathfrak h}$.
Recall that there exists
$H^j_{{\mathfrak k}} \in
({\mathfrak h})_{{\mathfrak k}_j} \ominus ({\mathfrak h} \cap {\mathfrak k}_j)$
such that
$H^j_{\mathfrak k} + H^j \in {\mathfrak h}$.
The next lemma shows that ${\mathfrak h}$ fails to be compatible with the Langlands
decomposition only for the abelian component ${\mathfrak a}_j$.

\begin{lm}
\label{lm:structure-h}
We have
\begin{itemize}
\item[(1)]
$\tau({\mathfrak h}) \cap {\mathfrak n}_j = {\mathfrak h} \cap {\mathfrak n}_j$,
\item[(2)]
${\mathfrak h} = ({\mathfrak h} \cap {\mathfrak m}_j)
\oplus \R (H^j_{\mathfrak k} + H^j) \oplus
({\mathfrak h} \cap {\mathfrak n}_j)$.
\end{itemize}
\end{lm}

\begin{proof}
Let $S$ be a connected solvable subgroup of $H$ which
acts transitively on the singular orbit $W$
(for the existence see e.g.\ Proposition 3.1 in \cite{BT3}),
and denote by ${\mathfrak s}$ the Lie algebra of $S$.
First of all, we show that
\begin{eqnarray}
\label{eq:structure-h-claim1}
\tau({\mathfrak s}) \cap {\mathfrak n}_j = {\mathfrak s} \cap {\mathfrak n}_j .
\end{eqnarray}
It is easy to see that
${\mathfrak s} \cap {\mathfrak n}_j
\subset \tau({\mathfrak s}) \cap {\mathfrak n}_j$.
We now choose
$Y_{\mathfrak n} \in (\tau({\mathfrak s}) \cap {\mathfrak n}_j)
\ominus ({\mathfrak s} \cap {\mathfrak n}_j)$,
and show that $Y_{\mathfrak n} = 0$.
Since $Y_{\mathfrak n} \in \tau({\mathfrak s})$,
there exists $Y_{\mathfrak k} \in {\mathfrak k}_j$ such that
$Y_{\mathfrak k} + Y_{\mathfrak n} \in {\mathfrak s}$.

We now claim that
\begin{eqnarray}
\label{eq:structure-h-claim2}
[ H^j_{\mathfrak k} , Y_{\mathfrak k} ] = 0 .
\end{eqnarray}
To show (\ref{eq:structure-h-claim2})
we define the solvable subalgebra
${\mathfrak s}^\prime := {\mathfrak s} \cap
({\mathfrak k}_j \oplus {\mathfrak a}_j \oplus {\mathfrak n}_j)$.
Let
$\pi_{\mathfrak k} :
{\mathfrak k}_j \oplus {\mathfrak a}_j \oplus {\mathfrak n}_j \to {\mathfrak k}_j$
be the canonical projection.
Since
${\mathfrak k}_j \subset {\mathfrak m}_j$ normalizes
${\mathfrak a}_j \oplus {\mathfrak n}_j$,
the map $\pi_{\mathfrak k}$ is a Lie algebra homomorphism,
and therefore $\pi_{\mathfrak k}({\mathfrak s}^\prime)$ is a solvable subalgebra of
${\mathfrak k}_j$.
Since every solvable subalgebra of a compact Lie algebra is abelian,
we conclude that $\pi_{\mathfrak k}({\mathfrak s}^\prime)$ is an
abelian subalgebra of ${\mathfrak k}_j$.
By construction, we have
$Y_{\mathfrak k} + Y_{\mathfrak n} , H^j_{\mathfrak k} + H^j
\in {\mathfrak s}^\prime$
and hence
$H^j_{\mathfrak k},Y_{\mathfrak k} \in \pi_{\mathfrak k}({\mathfrak s}^\prime)$.
This proves (\ref{eq:structure-h-claim2}).

Our next claim is
\begin{eqnarray}
\label{eq:structure-h-claim3}
[ H^j_{\mathfrak k} + H^j , Y_{\mathfrak n} ]
\in {\mathfrak s} \cap {\mathfrak n}_j .
\end{eqnarray}
Recall that
$H^j_{\mathfrak k} + H^j , Y_{\mathfrak k} + Y_{\mathfrak n} \in {\mathfrak s}$.
Hence we have
$$
{\mathfrak s} \ni
[ H^j_{\mathfrak k} + H^j , Y_{\mathfrak k} + Y_{\mathfrak n} ]
=
[ H^j_{\mathfrak k} , Y_{\mathfrak k} ]
+ [ H^j , Y_{\mathfrak k} ]
+ [ H^j_{\mathfrak k} + H^j , Y_{\mathfrak n} ]
= [ H^j_{\mathfrak k} + H^j , Y_{\mathfrak n} ] .
$$
Note that the last equality follows from
(\ref{eq:structure-h-claim2}) and
$[ {\mathfrak k}_j , {\mathfrak a}_j ] = 0$.
As
$H^j_{\mathfrak k} + H^j \in {\mathfrak k}_j \oplus {\mathfrak a}_j
\subset {\mathfrak l}_j$ and
${\mathfrak l}_j$ normalizes ${\mathfrak n}_j$,
we also have
$[H^j_{\mathfrak k} + H^j , Y_{\mathfrak n}] \in {\mathfrak n}_j$.
Thus (\ref{eq:structure-h-claim3}) has been proved.

Recall that
$Y_{\mathfrak n} \in (\tau({\mathfrak s}) \cap {\mathfrak n}_j)
\ominus ({\mathfrak s} \cap {\mathfrak n}_j)$.
From (\ref{eq:structure-h-claim3}) and the
skewsymmetry of $\ad(H^j_{{\mathfrak k}})$,
we have
$$
0 = \langle [ H^j_{\mathfrak k} + H^j , Y_{\mathfrak n} ] , Y_{\mathfrak n} \rangle
= \langle [ H^j , Y_{\mathfrak n} ] , Y_{\mathfrak n} \rangle .
$$
Recall that $H^j$ determines the gradation
${\mathfrak n}_j = \bigoplus_{\nu=1}^{m_j} {\mathfrak n}_j^\nu$,
and we therefore can write
$Y_{\mathfrak n} = \sum_{\nu=1}^{m_j} Y_{\mathfrak n}^\nu$
with $Y_{\mathfrak n}^\nu \in {\mathfrak n}_j^\nu$.
Hence we have
$$
0 = \langle [H^j , Y_{\mathfrak n} ] , Y_{\mathfrak n} \rangle
= \sum_{\nu=1}^{m_j} \nu
\langle Y_{\mathfrak n}^\nu , Y_{\mathfrak n}^\nu \rangle,
$$
which implies $Y_{\mathfrak n}^\nu = 0$ for all $\nu \in \{1,\ldots,m_j\}$.
We thus conclude that $Y_{\mathfrak n} = 0$,
and therefore (\ref{eq:structure-h-claim1}) has been proved.

We now prove statement (1) of the lemma.
It is easy to see
${\mathfrak h} \cap {\mathfrak n}_j \subset \tau({\mathfrak h}) \cap {\mathfrak n}_j$.
To show the converse,
note that $\tau({\mathfrak h}) = \tau({\mathfrak s})$
since $W = S \cdot o = H \cdot o$.
Hence, by (\ref{eq:structure-h-claim1}), we have
$$
\tau({\mathfrak h}) \cap {\mathfrak n}_j
= \tau({\mathfrak s}) \cap {\mathfrak n}_j
= {\mathfrak s} \cap {\mathfrak n}_j \subset {\mathfrak s} \subset {\mathfrak h} .
$$
This proves
$\tau({\mathfrak h}) \cap {\mathfrak n}_j \subset {\mathfrak h} \cap {\mathfrak n}_j$
and hence (1) holds.

We now prove statement (2) of the lemma.
It is easy to see ``$\supset$'' of (2).
To show ``$\subset$'', we choose $X \in {\mathfrak h}$
and write
$$
X = X_{\mathfrak m} + X_{\mathfrak a} + X_{\mathfrak n}
$$
according to the Langlands decomposition
${\mathfrak h} \subset {\mathfrak q}_j =
{\mathfrak m}_j \oplus {\mathfrak a}_j \oplus {\mathfrak n}_j$.
By definition of ${\mathfrak a}_j$
we can write $X_{\mathfrak a} = c H^j$ with some $c \in \R$,
and hence we can write $X$ as
$$
X = (X_{\mathfrak m} - c H^j_{\mathfrak k})
+ c (H^j_{\mathfrak k} + H^j) + X_{\mathfrak n} .
$$
One can easily see that
$$
X_{\mathfrak m} - c H^j_{\mathfrak k} \in {\mathfrak m}_j , \
c (H^j_{\mathfrak k} + H^j) \in \R (H^j_{\mathfrak k} + H^j) , \
X_{\mathfrak n} \in \tau({\mathfrak h}) \cap {\mathfrak n}_j .
$$
Note that
$X_{\mathfrak n} \in \tau({\mathfrak h}) \cap {\mathfrak n}_j
= {\mathfrak h} \cap {\mathfrak n}_j$
from (1).
Hence we have
$$
{\mathfrak h} \ni X - c (H^j_{\mathfrak k} + H^j) - X_{\mathfrak n}
= X_{\mathfrak m} - c H^j_{\mathfrak k} .
$$
This yields
$X_{\mathfrak m} - c H^j_{\mathfrak k} \in {\mathfrak h} \cap {\mathfrak m}_j$
and proves $X \in ({\mathfrak h} \cap {\mathfrak m}_j)
\oplus \R (H^j_{\mathfrak k} + H^j) \oplus
({\mathfrak h} \cap {\mathfrak n}_j)$.
This finishes the proof of (2).
\end{proof}

We next show that ${\mathfrak h}$ can be replaced by a simpler subalgebra with an orbit equivalent action.

\begin{lm}
\label{lm:h-orbit-eq}
The action of $H$ is orbit equivalent to the action of the connected Lie subgroup
$H^\prime$ of $Q_j$ with Lie algebra
${\mathfrak h}^\prime := ({\mathfrak h} \cap {\mathfrak m}_j)
\oplus {\mathfrak a}_j \oplus ({\mathfrak h} \cap {\mathfrak n}_j)$.
\end{lm}

\begin{proof}
From Proposition \ref{normalspace}
we know that
$\nu_oW \subset {\mathfrak b}_j$ or $\nu_oW \subset {\mathfrak n}_j^1$,
and thus we have
\begin{eqnarray}
\label{eq:tau-h-cap-n}
{\mathfrak h} \cap {\mathfrak n}_j
= {\mathfrak n}_j
\ \ \mbox{or} \ \
{\mathfrak h} \cap {\mathfrak n}_j
= {\mathfrak n}_j \ominus \nu_o W
= ({\mathfrak n}^1_j \ominus \nu_o W) \oplus
(\bigoplus_{\nu > 1} {\mathfrak n}^{\nu}_j) .
\end{eqnarray}

First of all, we show that ${\mathfrak h}^\prime$ is a subalgebra.
It is easy to see that
$({\mathfrak h} \cap {\mathfrak m}_j) \oplus {\mathfrak a}_j$ is a subalgebra,
and from (\ref{eq:tau-h-cap-n})
we see that
${\mathfrak a}_j \oplus ({\mathfrak h} \cap {\mathfrak n}_j)$
is a subalgebra.
Since ${\mathfrak m}_j$ normalizes ${\mathfrak n}_j$, we also have
$[ {\mathfrak h} \cap {\mathfrak m}_j , {\mathfrak h} \cap {\mathfrak n}_j ]
\subset {\mathfrak h} \cap {\mathfrak n}_j $.
Altogether this implies that ${\mathfrak h}^\prime$ is a subalgebra.

Next we prove that $\R H^j_{\mathfrak k} \oplus {\mathfrak h}^\prime$
is a subalgebra.
Since ${\mathfrak h}^\prime$ is a subalgebra,
it is enough to show that
\begin{itemize}
\item[(i)]
$[ H^j_{\mathfrak k} , {\mathfrak h} \cap {\mathfrak m}_j ]
\subset {\mathfrak h} \cap {\mathfrak m}_j$, \
(ii)
$[ H^j_{\mathfrak k} , {\mathfrak a}_j ] = 0$, \
(iii)
$[ H^j_{\mathfrak k} , {\mathfrak h}\cap {\mathfrak n}_j ]
\subset {\mathfrak h} \cap {\mathfrak n}_j$.
\end{itemize}
Let $X \in {\mathfrak h} \cap {\mathfrak m}_j$.
Since $H^j_{\mathfrak k} \in {\mathfrak k}_j \subset {\mathfrak m}_j$ and
${\mathfrak m}_j$ is a subalgebra,
we have
$[ H^j_{\mathfrak k} , X ] \in {\mathfrak m}_j$.
Furthermore, since $H^j$ centralizes ${\mathfrak m}_j$, we have
$
[ H^j_{\mathfrak k} , X ] =
[ H^j_{\mathfrak k} + H^j , X ] \in {\mathfrak h} $.
Altogether this gives
$[ H^j_{\mathfrak k} , X ] \in {\mathfrak h} \cap {\mathfrak m}_j$,
which implies (i).
The claim (ii) is easy to verify.
If $ {\mathfrak h} \cap {\mathfrak n}_j = {\mathfrak n}_j$,
then it is easy to see that (iii) holds.
If $ {\mathfrak h} \cap {\mathfrak n}_j
= {\mathfrak n}_j \ominus \nu_o W$,
then $\ad(H^j_{{\mathfrak k}})$ normalizes $\nu_o W$ by Lemma \ref{lm:f^2},
and hence normalizes ${\mathfrak n}_j \ominus \nu_o W$ by
skewsymmetry of $\ad(H^j_{{\mathfrak k}})$.
This finishes the proof of (iii).

We now consider the three subalgebras
${\mathfrak h}$,
$\R H^j_{\mathfrak k} \oplus {\mathfrak h}^\prime$ and
${\mathfrak h}^\prime$.
By construction,  we have
${\mathfrak h} ,  {\mathfrak h}^\prime
\subset \R H^j_{\mathfrak k} \oplus {\mathfrak h}^\prime$.
Denote by $H''$ the connected Lie subgroup of $Q_j$ with Lie algebra
$\R H^j_{\mathfrak k} \oplus {\mathfrak h}^\prime$.
We will now prove that the actions of $H$, $H^\prime$ and $H''$ are orbit equivalent
to each other.

We first show that the actions of $H$ and $H''$ are orbit equivalent.
Since $H \subset H''$,
we obviously have $W = H \cdot o \subset H'' \cdot o$.
However, since
$T_o W = \tau({\mathfrak h}) = \tau(\R H^j_{\mathfrak k} \oplus {\mathfrak h}^\prime)$
and both orbits are connected and complete,
we conclude that $W = H \cdot o = H'' \cdot o$.
By assumption the action of $H$ is of cohomogeneity one,
and therefore the action of $H''$ must be of cohomogeneity one as well.
This implies that the actions of $H$ and $H''$ are orbit equivalent.

We next show that the actions of $H^\prime$ and $H''$ are orbit equivalent.
Since $H^\prime \subset H''$ and
$\tau({\mathfrak h}') = \tau(\R H^j_{\mathfrak k} \oplus {\mathfrak h}^\prime)
= T_o W$,
we conclude that $H^\prime \cdot o = H'' \cdot o = W$.
By construction, we have
${\mathfrak h}' \cap {\mathfrak k}_j = {\mathfrak h} \cap {\mathfrak k}_j $,
which implies that the slice representations of $H^\prime$ and $H$ at $o$ are the same.
Since the action of $H$ is of cohomogeneity one by assumption,
the action of $H'$ is of cohomogeneity one as well.
Thus, since both  actions are of cohomogeneity one
and have the same singular orbit,
we conclude that these actions are orbit equivalent.

We thus have proved that the actions of $H$ and $H^\prime$ are orbit equivalent.
\end{proof}

According to Proposition \ref{normalspace},
the normal space of the singular orbit is either tangent to the
(totally geodesic) semisimple part or to the nilpotent part of the horospherical
decomposition of $M$ induced by $\Phi_j$. We now distinguish these two cases.

\begin{prop}
If $\nu_oW \subset {\mathfrak b}_j$, then the action of $H$ on $M$
is orbit equivalent to the canonical extension of
a cohomogeneity one action on the boundary component $B_j$ of $M$.
\end{prop}

\begin{proof}
Assume that $\nu_oW \subset {\mathfrak b}_j$. According to Lemma \ref{lm:h-orbit-eq} we can assume that
$$
{\mathfrak h} = ({\mathfrak h} \cap {\mathfrak m}_j)
\oplus {\mathfrak a}_j \oplus {\mathfrak n}_j .
$$
Note that ${\mathfrak m}_j$ is reductive and we have the Lie algebra
direct sum decomposition
${\mathfrak m}_j = {\mathfrak z}_j \oplus {\mathfrak g}_j$,
where ${\mathfrak z}_j$ is the center of ${\mathfrak m}_j$.
Therefore the canonical projection
$\pi_{\mathfrak g} : {\mathfrak m}_j \to {\mathfrak g}_j$
with respect to this decomposition is a Lie algebra homomorphism and
${\mathfrak h}^\prime := \pi_{\mathfrak g}({\mathfrak h} \cap {\mathfrak m}_j)$
is a subalgebra of ${\mathfrak g}_j$.
Let $H^\prime$ be the connected subgroup of $G_j$ with Lie algebra ${\mathfrak h}^\prime$.
We claim that $H^\prime$ acts on $B_j = M_j/K_j = G_j/(G_j \cap K_j)$
with cohomogeneity one and the canonical extension of this action
to $M$ is orbit equivalent to the action of $H$ on $M$.

We first prove that $H^\prime$ acts on $B_j$ with cohomogeneity one.
For simplicity we will identify the subalgebras and the corresponding connected Lie subgroups.
At first we consider the action of ${\mathfrak h} \cap {\mathfrak m}_j$ on $B_j$.
The slice representation of this action is the action of
${\mathfrak h} \cap {\mathfrak k}_j$ on $\nu_o W$,
which coincides with the slice representation of the action of $H$ on $M$.
Therefore, ${\mathfrak h} \cap {\mathfrak k}_j$ acts transitively on
the unit sphere in $\nu_o W$,
and hence the action of ${\mathfrak h} \cap {\mathfrak m}_j$ on $B_j$
is of cohomogeneity one.
Next we consider
$$
{\mathfrak h} \cap {\mathfrak m}_j
\subset {\mathfrak h}^\prime \oplus {\mathfrak z}_j .
$$
Since $\tau({\mathfrak h} \cap {\mathfrak m}_j)
= \tau({\mathfrak h}^\prime \oplus {\mathfrak z}_j)$,
the action of ${\mathfrak h}^\prime \oplus {\mathfrak z}_j$ on $B_j$
is also of cohomogeneity one.
Finally we consider
$$
{\mathfrak h}^\prime \oplus {\mathfrak z}_j \supset {\mathfrak h}^\prime .
$$
Since $\tau({\mathfrak h} \cap {\mathfrak m}_j) = \tau({\mathfrak h}^\prime)$,
their orbits through $o$ coincide.
Furthermore, since ${\mathfrak z}_j$ acts trivially on $\nu_o W$,
the slice representations of these actions are equivalent.
Therefore we conclude that the action of $H^\prime$ on $B_j$ is of
cohomogeneity one.

We now consider the canonical extension $H_j^\Lambda$ of $H^\prime$ to $M$.
By definition, we have
$$
{\mathfrak h}_j^\Lambda =
{\mathfrak h}^\prime \oplus {\mathfrak a}_j \oplus {\mathfrak n}_j .
$$
By a similar argument as above, one can show that
the following three actions are orbit equivalent:
$$
{\mathfrak h}_j^\Lambda \subset
{\mathfrak h}_j^\Lambda \oplus {\mathfrak z}_j \supset {\mathfrak h} .
$$
Therefore,
the action of $H$ is orbit equivalent to the action of the canonical extension
$H_j^\Lambda$ of $H^\prime$.
\end{proof}

We now turn our attention to cohomogeneity one actions with $\nu_oW \subset {\mathfrak n}_j^1$.

\begin{prop}
If $\nu_oW \subset {\mathfrak n}_j^1$,
then the action of $H$ on $M$
is orbit equivalent to the action of $H_{j,{\mathfrak v}}$ for some
subspace ${\mathfrak v} \subset {\mathfrak n}_j^1$.
\end{prop}

\begin{proof}
Assume that $\nu_oW \subset {\mathfrak n}_j^1$.
By means of Lemma \ref{lm:h-orbit-eq} we can assume that
$$
{\mathfrak h} = ({\mathfrak h} \cap {\mathfrak m}_j)
\oplus {\mathfrak a}_j \oplus ({\mathfrak n}_j \ominus \nu_o W) .
$$
Let ${\mathfrak v} := \nu_o W$ and recall that
$$
{\mathfrak h}_{j,{\mathfrak v}}
= N_{{\mathfrak l}_j} ({\mathfrak n}_{j,{\mathfrak v}})
\oplus {\mathfrak n}_{j,{\mathfrak v}}
= N_{{\mathfrak l}_j} ({\mathfrak n}_{j,{\mathfrak v}})
\oplus ({\mathfrak n}_j \ominus \nu_o W) .
$$
We have
${\mathfrak h} \subset {\mathfrak h}_{j,{\mathfrak v}}$
since
${\mathfrak h}$ is a subalgebra and hence
$({\mathfrak h} \cap {\mathfrak m}_j) \oplus {\mathfrak a}_j$
normalizes ${\mathfrak n}_{j,{\mathfrak v}}$.
One can also see that
$\tau({\mathfrak h}) = \tau({\mathfrak h}_{j,{\mathfrak v}})$, and
therefore
$H \cdot o = H_{j,{\mathfrak v}} \cdot o$.
Since the actions of $H$ and $H_{j,{\mathfrak v}}$ are of cohomogeneity one,
these actions are orbit equivalent.
\end{proof}

From the previous two propositions we obtain the main result of this section.

\begin{thm} \label{parabolicmaintheorem}
Let $M$ be a connected Riemannian symmetric space of noncompact type,
and let $H$ be a connected subgroup of $G$ which acts on $M$ with cohomogeneity one
and has a singular orbit. Assume that $H$ is
contained in a proper maximal parabolic subgroup $Q_j$ of $G$.
Then the action of $H$ on $M$ is orbit equivalent to
\begin{itemize}
\item[(i)] a cohomogeneity one action on $M$ obtained by canonical extension of a cohomogeneity
one action on the maximal boundary component $B_j$, or
\item[(ii)] a cohomogeneity one action on $M$ given by $H_{j,{\mathfrak v}}$ for some
subspace ${\mathfrak v} \subset {\mathfrak n}_j^1$.
\end{itemize}
\end{thm}

We emphasize that in Theorem \ref{parabolicmaintheorem} the symmetric space $M$ can be reducible.
As a consequence of this result we also see that a non-totally geodesic singular orbit of a cohomogeneity one
action on $M$ which is not a canonical extension contains a maximal boundary component of $M$.

This finishes the proof of Theorem \ref{maintheorem}.

\section{Some explicit classifications}

In this section we present explicit classifications of cohomogeneity
one actions (up to orbit equivalence) for some symmetric spaces of rank two.
We have chosen symmetric spaces for which the Lie algebra ${\mathfrak g}$ of the
isometry group is a split real form of its complexification ${\mathfrak g}^{\mathbb C}$.
In order to describe these
we recall briefly the classification of cohomogeneity one actions on a real hyperbolic space
${\mathbb R}H^n$ (see \cite{BT3} for further details).

\begin{thm} \label{hyperbolic_C1}
Every cohomogeneity one action on the real hyperbolic space ${\mathbb R}H^n = SO^o_{1,n}/SO_n$
is orbit equivalent to one of the following actions:
\begin{itemize}
\item[(1)] The action of $SO^o_{1,k} \times SO_{n-k} \subset SO^o_{1,n}$ for some $k \in
\{0,\ldots,n-1\}$. For $k < n-1$ the action has exactly one singular orbit, namely a
totally geodesic ${\mathbb R}H^k \subset {\mathbb R}H^n$. For $k = n-1$ the orbits form
a foliation, and one of the orbits is a totally geodesic ${\mathbb R}H^{n-1}$.
\item[(ii)] The action of the nilpotent subgroup $N$ in an Iwasawa decomposition
$SO^o_{1,n} = SO_nAN$. The orbits form a foliation on ${\mathbb R}H^n$ by horospheres.
\end{itemize}
\end{thm}

Since $SO^o_{1,n}$ acts transitively on ${\mathbb R}H^n$ and the isotropy group at a point
is isomorphic to $SO_n$, it can easily be seen that orbit equivalence can always be achieved
 by an isometry in $SO^o_{1,n}$.

\subsection{The symmetric space $M = SL_3({\mathbb R})/SO_3$.}
The symmetric space $M = SL_3({\mathbb R})/SO_3$ has rank $2$
and dimension $5$.
The root system is of type $(A_2)$ and all multiplicities are equal to $1$.
The positive roots are $\alpha_1,\alpha_2,\alpha_1+\alpha_2$ and the nilpotent subalgebra
${\mathfrak n}$ of ${\mathfrak g} = {\mathfrak s}{\mathfrak l}_3({\mathbb R})$ is given by
$$
{\mathfrak n} = {\mathfrak g}_{\alpha_1} \oplus {\mathfrak g}_{\alpha_2} \oplus
{\mathfrak g}_{\alpha_1+\alpha_2}.
$$
The maximal abelian subalgebra ${\mathfrak a}$ has dimension $2$ and is spanned
by the two root vectors $H_{\alpha_1}$ and $H_{\alpha_2}$.
The Chevalley decomposition ${\mathfrak q}_2 = {\mathfrak l}_2  \oplus {\mathfrak n}_2$
of the parabolic subalgebra ${\mathfrak q}_2$ corresponding to $\Phi_2 = \{\alpha_1\}$ is given by
$${\mathfrak l}_2 = {\mathfrak g}_{-\alpha_1} \oplus {\mathfrak g}_0 \oplus
{\mathfrak g}_{\alpha_1} \cong {\mathfrak s}{\mathfrak l}_2({\mathbb R}) \oplus {\mathbb R}
\ \ {\rm and}\ \ {\mathfrak n}_2 = {\mathfrak g}_{\alpha_2} \oplus {\mathfrak g}_{\alpha_1+\alpha_2}.$$
The orbit $F_2 = L_2 \cdot o$ is isometric to the Riemannian product
${\mathbb R}H^2 \times {\mathbb E}$,
and the corresponding boundary component $B_2$ is the real hyperbolic plane ${\mathbb R}H^2$.

\begin{thm} \label{A21}
Each cohomogeneity one action on
$M = SL_3({\mathbb R})/SO_3$
is orbit equivalent to one of the following cohomogeneity one actions on
$M$:
\begin{itemize}
\item[(1)] The action of the subgroup $H_\ell$ of $SL_3({\mathbb R})$ with Lie algebra
$${\mathfrak h}_\ell = ({\mathfrak a} \ominus \ell) \oplus {\mathfrak n},$$
where $\ell$ is a one-dimensional linear subspace of ${\mathfrak a}$. The orbits form
a Riemannian foliation on $M$ and all orbits are isometrically congruent to each other.
\item[(2)] The action of the subgroup $H_1$ of $SL_3({\mathbb R})$ with Lie algebra
$${\mathfrak h}_1 = ({\mathfrak a} \oplus {\mathfrak n}) \ominus {\mathfrak g}_{\alpha_1}.$$
The orbits form a Riemannian foliation on $M$ and there is exactly one
minimal orbit $H_1 \cdot o$.
\item[(3)] The action of $SL_2({\mathbb R}) \times {\mathbb R}^+ \subset SL_3({\mathbb R})$
with Lie algebra
$${\mathfrak l}_2 = {\mathfrak g}_{-\alpha_1} \oplus {\mathfrak g}_0 \oplus {\mathfrak g}_{\alpha_1} \cong
{\mathfrak s}{\mathfrak l}_2({\mathbb R}) \oplus {\mathbb R}.$$
This action has a totally geodesic singular orbit isometric to ${\mathbb R}H^2 \times {\mathbb E}$.
\item[(4)] The action of the connected subgroup $H$ of $SL_3({\mathbb R})$ with Lie algebra
$${\mathfrak h} = {\mathfrak k}_{\alpha_1} \oplus ({\mathfrak a} \ominus
{\mathbb R}H_{\alpha_1}) \oplus ({\mathfrak n} \ominus {\mathfrak g}_{\alpha_1}) ,$$
where ${\mathfrak k}_{\alpha_1} \cong {\mathfrak s}{\mathfrak o}_2$ is the Lie algebra of the
isotropy group of the isometry group of the boundary component $B_2 = {\mathbb R}H^2$.
This action has a minimal ${\mathbb R}H^3 \subset M$ as a singular orbit and can be
constructed by canonical extension of the
cohomogeneity one action on $B_2 = {\mathbb R}H^2$ of $M$ which has
a single point as an orbit.
\end{itemize}
\end{thm}

\begin{proof}
For the classification we have to consider the different cases in Theorem \ref{maintheorem}.
If the orbits form a Riemannian foliation, we obtain the actions described in (1) and (2).
The action described in (3) is the only one corresponding to case (2)(i) in Theorem \ref{maintheorem}
according to \cite{BT2}.
We now consider an action as described in Theorem \ref{maintheorem} (2)(ii).
The symmetric space $M$ has, up to isometric congruence, only one boundary component with
rank one, namely the real hyperbolic plane $B_2 = {\mathbb R}H^2$. There is, up to orbit equivalence,
exactly one cohomogeneity one action on ${\mathbb R}H^2$ with a  singular orbit,
namely the action on ${\mathbb R}H^2$ by the isotropy group $K_{\alpha_1} \cong SO_2$. The
canonical extension of this action leads to the action described in (4).
It remains to investigate case (b) in (2)(ii) of Theorem \ref{maintheorem}. For $\Phi_2 = \{\alpha_1\}$
we have
${\mathfrak n}_2^1 = {\mathfrak n}_2 = {\mathfrak g}_{\alpha_2} \oplus
{\mathfrak g}_{\alpha_1+\alpha_2}$,
and therefore ${\mathfrak v} = {\mathfrak n}_2^1$ for dimension reasons. It is easy to see that
$N^o_{K_2}({\mathfrak v}) = K_2 \cong SO_2$ acts transitively on the unit sphere in ${\mathfrak v}$. Moreover,
the normalizer of ${\mathfrak n}_2 \ominus {\mathfrak v} = \{0\}$ in $L_2$ is clearly
$L_2$, which acts transitively on $F_2 = {\mathbb R}H^2 \times {\mathbb E}$. The construction method
in (2)(ii)(b) therefore leads to the cohomogeneity one action on $M$ by $L_2$, which is the action
described in (3). The case $\Phi_1 = \{\alpha_2\}$ does not lead to anything new because of the Dynkin diagram symmetry.
\end{proof}

\subsection{The symmetric space $SO^o_{2,3}/SO_2SO_3 = G_2^*({\mathbb R}^5)$.}
The symmetric space $M = SO^o_{2,3}/SO_2SO_3$ has rank $2$ and dimension $6$.
The root system is of type $(B_2)$, the Dynkin diagram is
$$
\xy
\POS
(40,0) *\cir<2pt>{}="d",
(50,0) *\cir<2pt>{}="e",
(40,-5) *{\alpha_1},
(50,-5) *{\alpha_2},
\ar @2{->} "d";"e",
\endxy
$$
and all multiplicities are equal to $1$.
The maximal abelian subalgebra ${\mathfrak a}$ has dimension $2$ and is spanned
by the two root vectors $H_{\alpha_1}$ and $H_{\alpha_2}$. Both boundary components $B_1$ and $B_2$
are isometric to a real hyperbolic plane ${\mathbb R}H^2$. However $B_1$ and $B_2$ are not isometrically
congruent to each other.

\begin{thm} \label{B21}
Each cohomogeneity one action on
$M = SO^o_{2,3}/SO_2SO_3$
is orbit equivalent to one of the following cohomogeneity one actions on
$M$:
\begin{itemize}
\item[(1)] The action of the subgroup $H_\ell$ of $SO^o_{2,3}$ with Lie algebra
$${\mathfrak h}_\ell = ({\mathfrak a} \ominus \ell) \oplus {\mathfrak n},$$
where $\ell$ is a one-dimensional linear subspace of ${\mathfrak a}$. The orbits form
a Riemannian foliation on $M$ and all orbits are isometrically congruent to each other.
\item[(2)] The action of the subgroup $H_i$, $i \in \{1,2\}$, of $SO^o_{2,3}$ with Lie algebra
$${\mathfrak h}_i = ({\mathfrak a} \oplus {\mathfrak n}) \ominus {\mathfrak g}_{\alpha_i}.$$
The orbits form a Riemannian foliation on $M$ and there is exactly one
minimal orbit $H_i \cdot o$.
\item[(3)] The action of $SO^o_{1,3} \subset SO^o_{2,3}$.
This action has a totally geodesic singular orbit isometric to the real hyperbolic space
${\mathbb R}H^3$.
\item[(4)] The action of $SO^o_{2,2} \subset SO^o_{2,3}$.
This action has a totally geodesic singular orbit isometric to the Riemannian product
${\mathbb R}H^2 \times {\mathbb R}H^2$.
\item[(5)] The action of the subgroup $H_1^\Lambda$ of $SO^o_{2,3}$ with Lie algebra
$${\mathfrak h}_1^\Lambda = {\mathfrak k}_{\alpha_2} \oplus ({\mathfrak a} \ominus {\mathbb R}H_{\alpha_2})
\oplus ({\mathfrak n} \ominus {\mathfrak g}_{\alpha_2}) ,$$
where ${\mathfrak k}_{\alpha_2} \cong {\mathfrak s}{\mathfrak o}_2$ is the Lie algebra of the
isotropy group of the isometry group of the boundary component $B_1 \cong {\mathbb R}H^2$.
The action of $H_1^\Lambda$ has a minimal real hyperbolic space ${\mathbb R}H^4 \subset M$ as
a singular orbit and can be
constructed by canonical extension of the
cohomogeneity one action on the boundary component $B_1 = {\mathbb R}H^2$ which has
a single point as an orbit.
\item[(6)] The action of the subgroup $H_2^\Lambda$ of $SO^o_{2,3}$ with Lie algebra
$${\mathfrak h}_2^\Lambda = {\mathfrak k}_{\alpha_1} \oplus ({\mathfrak a} \ominus {\mathbb R}H_{\alpha_1})
\oplus ({\mathfrak n} \ominus {\mathfrak g}_{\alpha_1}) ,$$
where ${\mathfrak k}_{\alpha_1} \cong {\mathfrak s}{\mathfrak o}_2$ is the Lie algebra of the
isotropy group of the isometry group of the boundary component $B_2 \cong {\mathbb R}H^2$.
The action of $H_2^\Lambda$ has a minimal complex hyperbolic plane ${\mathbb C}H^2 \subset M$
as a singular orbit and can be constructed by canonical extension of the
cohomogeneity one action on the boundary component $B_2 = {\mathbb R}H^2$ which has
a single point as an orbit.
\end{itemize}
\end{thm}

\begin{proof}
For the classification we have to consider the different cases in Theorem \ref{maintheorem}.
If the orbits form a Riemannian foliation, we obtain the actions  in (1) and (2).
The actions  in (3) and (4) are the only ones corresponding to case (2)(i) in Theorem \ref{maintheorem}
according to \cite{BT2}.
We now consider an action as described in Theorem \ref{maintheorem} (2)(ii).
The symmetric space $M$ has two maximal boundary components $B_1$ and $B_2$. Both $B_1$ and $B_2$
are isometric to ${\mathbb R}H^2$ with a suitable constant curvature metric,
but they are not isometrically congruent in $M$.
There is, up to orbit equivalence,
exactly one cohomogeneity one action on ${\mathbb R}H^2$ with a  singular orbit,
namely the action on ${\mathbb R}H^2$ by the isotropy group $SO_2$. The
canonical extension of this action leads to the actions  in (5) and (6).
It remains to investigate case (b) in (2)(ii) of Theorem \ref{maintheorem}.
We have to consider two possible choices of subsystems of $\Lambda = \{\alpha_1,\alpha_2\}$,
namely $\Phi_1 = \{\alpha_2\}$ and
$\Phi_2 = \{\alpha_1\}$.

In case of $\Phi_1$ we have
\begin{eqnarray*}
{\mathfrak n}_1^1 & = &
{\mathfrak g}_{\alpha_1} \oplus {\mathfrak g}_{\alpha_1 + \alpha_2}
\oplus {\mathfrak g}_{\alpha_1 + 2 \alpha_2} \cong {\mathbb R}^3,\\
{\mathfrak l}_1 & = & {\mathfrak g}_{-\alpha_2} \oplus {\mathfrak g}_0 \oplus {\mathfrak g}_{\alpha_2}
\cong {\mathfrak s}{\mathfrak l}_2({\mathbb R}) \oplus {\mathbb R},\\
{\mathfrak k}_1 & = & {\mathfrak k}_{\alpha_2} \cong {\mathfrak s}{\mathfrak o}_2.
\end{eqnarray*}
Since ${\mathfrak k}_1$ is one-dimensional,
${\mathfrak v}$ must be a $2$-dimensional linear subspace of ${\mathfrak n}_1^1$.
Let ${\mathfrak v}$ be a ${\mathfrak k}_1$-invariant
$2$-dimensional subspace of ${\mathfrak n}_1^1$.
In order to get a cohomogeneity one action, the normalizer $N_{L_1}({\mathfrak n}_1^1 \ominus {\mathfrak v})$
must act transitively on $F_1 = L_1 \cdot o \cong {\mathbb R}H^2 \times {\mathbb E}$. The only subgroups
of $SL_2({\mathbb R})$ acting transitively on ${\mathbb R}H^2$ are $SL_2({\mathbb R})$ itself and
the parabolic subgroups of $SL_2({\mathbb R})$. However,
${\mathfrak m}_1 \cong {\mathfrak s}{\mathfrak l}_2({\mathbb R})$ acts
irreducibly on ${\mathfrak n}_1^1$, and hence $N_{L_2}({\mathfrak n}_1^1 \ominus {\mathfrak v})$ cannot
be equal to $SL_2({\mathbb R})$. Since $K_1$ is compact and normalizes ${\mathfrak v}$, it also
normalizes ${\mathfrak n}_1^1 \ominus {\mathfrak v}$. If a parabolic subgroup of $SL_2({\mathbb R})$
would normalize ${\mathfrak n}_1^1 \ominus {\mathfrak v}$, then the entire group $SL_2({\mathbb R})$
would normalize ${\mathfrak n}_1^1 \ominus {\mathfrak v}$, which cannot happen. We thus conclude that
$N_{L_1}({\mathfrak n}_1^1 \ominus {\mathfrak v})$ cannot act transitively on $F_1$. This implies that there
is no cohomogeneity one action on $M$ which can be constructed from the choice of $\Phi_1$.

In case of $\Phi_2$ we have
\begin{eqnarray*}
{\mathfrak n}_2^1 & = &
{\mathfrak g}_{\alpha_2} \oplus {\mathfrak g}_{\alpha_1 + \alpha_2} \cong {\mathbb R}^2,\\
{\mathfrak n}_2^2 & = & {\mathfrak g}_{\alpha_1 + 2 \alpha_2} \cong {\mathbb R},\\
{\mathfrak l}_2 & = & {\mathfrak g}_{-\alpha_1} \oplus {\mathfrak g}_0 \oplus {\mathfrak g}_{\alpha_1}
\cong {\mathfrak s}{\mathfrak l}_2({\mathbb R}) \oplus {\mathbb R},\\
{\mathfrak k}_2 & = & {\mathfrak k}_{\alpha_1} \cong {\mathfrak s}{\mathfrak o}_1.
\end{eqnarray*}
The only possible choice for ${\mathfrak v}$ is therefore ${\mathfrak v} = {\mathfrak n}_2^1$.
The normalizer $N_{L_2}({\mathfrak n}_2^1 \ominus {\mathfrak v})$ is of course $L_2$, and
therefore we get a cohomogeneity one action on $M$. However, this action has $(L_2N_2^2) \cdot o
\cong {\mathbb R}H^2 \times {\mathbb R}H^2$ as a totally geodesic singular orbit, which we already
listed in (4).
\end{proof}

\subsection{The symmetric space $G_2^2/SO_4$.} The symmetric space $M = G_2^2/SO_4$ has rank $2$
and dimension $8$.
The root system is of type $(G_2)$, the Dynkin diagram is
$$
\xy
\POS (0,0) *\cir<2pt>{} ="a",
(10,0) *\cir<2pt>{}="b",
(0,-5) *{\alpha_1},
(10,-5) *{\alpha_2},
\ar @3{<-} "a";"b",
\endxy,
$$
and all multiplicities are equal to $1$.
The maximal abelian subalgebra ${\mathfrak a}$ has dimension two and is spanned
by the two root vectors $H_{\alpha_1}$ and $H_{\alpha_2}$. Both boundary components $B_1$ and $B_2$
are isometric to a real hyperbolic plane ${\mathbb R}H^2$. However, $B_1$ and $B_2$ are not isometrically
congruent to each other.

\begin{thm} \label{G21}
Each cohomogeneity one action on $M = G_2^2/SO_4$
is orbit equivalent to one of the following cohomogeneity one actions on
$M$:
\begin{itemize}
\item[(1)] The action of the subgroup $H_\ell$ of $G_2^2$ with Lie algebra
$${\mathfrak h}_\ell = ({\mathfrak a} \ominus \ell) \oplus {\mathfrak n},$$
where $\ell$ is a one-dimensional linear subspace of ${\mathfrak a}$. The orbits form
a Riemannian foliation on $M$ and all orbits are isometrically congruent to each other.
\item[(2)] The action of the subgroup $H_i$, $i \in \{1,2\}$, of $G_2^2$ with Lie algebra
$${\mathfrak h}_i = ({\mathfrak a} \oplus {\mathfrak n}) \ominus {\mathfrak g}_{\alpha_i}.$$
The orbits form a Riemannian foliation on $M$ and there is exactly one
minimal orbit $H_i \cdot o$.
\item[(3)] The action of $SU_{1,2} \subset G_2^2$.
This action has a totally geodesic singular orbit isometric to the complex hyperbolic plane
${\mathbb C}H^2$.
\item[(4)] The action of $SL_3({\mathbb R}) \subset G_2^2$.
This action has a totally geodesic singular orbit isometric to the symmetric space
$SL_3({\mathbb R})/SO_3$.
\item[(5)] The action of the subgroup $H_1^\Lambda$ of $G_2^2$ with Lie algebra
$${\mathfrak h}_1^\Lambda = {\mathfrak k}_{\alpha_2} \oplus ({\mathfrak a} \ominus {\mathbb R}H_{\alpha_2})
\oplus ({\mathfrak n} \ominus {\mathfrak g}_{\alpha_2}) ,$$
where ${\mathfrak k}_{\alpha_2} \cong {\mathfrak s}{\mathfrak o}_2$ is the Lie algebra of the
isotropy group of the isometry group of the boundary component $B_1 \cong {\mathbb R}H^2$.
This action has a $6$-dimensional minimal singular orbit and can be
constructed by canonical extension of the
cohomogeneity one action on the boundary component $B_1 = {\mathbb R}H^2$ which has
a single point as an orbit.
\item[(6)] The action of the subgroup $H_2^\Lambda$ of $SO^o_{2,3}$ with Lie algebra
$${\mathfrak h}_2^\Lambda = {\mathfrak k}_{\alpha_1} \oplus ({\mathfrak a} \ominus {\mathbb R}H_{\alpha_1})
\oplus ({\mathfrak n} \ominus {\mathfrak g}_{\alpha_1}) ,$$
where ${\mathfrak k}_{\alpha_1} \cong {\mathfrak s}{\mathfrak o}_2$ is the Lie algebra of the
isotropy group of the isometry group of the boundary component $B_2 \cong {\mathbb R}H^2$.
This action has a minimal complex hyperbolic space ${\mathbb C}H^3 \subset M$
as a singular orbit and can be
constructed by canonical extension of the
cohomogeneity one action on the boundary component $B_2 = {\mathbb R}H^2$ which has
a single point as an orbit.
\item[(7)] The action of the subgroup $H_{1,{\mathfrak v}}$ of $G_2^2$ with ${\mathfrak v} =
{\mathfrak g}_{\alpha_1} \oplus {\mathfrak g}_{\alpha_1 + \alpha_2}$ and Lie algebra
$$
{\mathfrak h}_{1,{\mathfrak v}} = {\mathfrak g}_{-\alpha_2} \oplus
{\mathfrak g}_0 \oplus {\mathfrak g}_{\alpha_2}
\oplus {\mathfrak g}_{2\alpha_1 + \alpha_2} \oplus
{\mathfrak g}_{3\alpha_1+\alpha_2} \oplus {\mathfrak g}_{3\alpha_1
+ 2\alpha_2}.
$$
This action has a $6$-dimensional minimal singular orbit.
\end{itemize}
\end{thm}

\begin{proof}
The argumentation for cases (1) to (6) is analogous to the one given in the proof of Theorem \ref{B21}.
We now consider the two possible choices of subsystems of $\Lambda = \{\alpha_1,\alpha_2\}$,
namely $\Phi_1 = \{\alpha_2\}$ and $\Phi_2 = \{\alpha_1\}$.

In case of $\Phi_1$ we have ${\mathfrak n}_1^1 =
{\mathfrak g}_{\alpha_1} \oplus {\mathfrak g}_{\alpha_1 + \alpha_2} \cong {\mathbb R}^2$.
The only possible choice for ${\mathfrak v}$ is therefore ${\mathfrak v} = {\mathfrak n}_1^1$.
This was discussed in detail in subsection \ref{method2},
where we showed that this leads to the cohomogeneity
one action described in (7). This action cannot be orbit equivalent to the one in (5) or (6),
as it contains a maximal flat of $M$ , whereas the two singular orbits in (5) and (6) do not contain a maximal
flat of $M$.

Finally, we consider $\Phi_2$. In this case we have
\begin{eqnarray*}
{\mathfrak n}_2^1 & = & {\mathfrak g}_{\alpha_2}
\oplus {\mathfrak g}_{\alpha_1 + \alpha_2} \oplus
{\mathfrak g}_{2\alpha_1+\alpha_2} \oplus {\mathfrak g}_{3\alpha_1
+ \alpha_2} \cong {\mathbb R}^4\\
{\mathfrak n}_2^2 & = & {\mathfrak g}_{3\alpha_1 + 2\alpha_2}\\
{\mathfrak l}_2 & = & {\mathfrak g}_{-\alpha_1} \oplus {\mathfrak g}_0 \oplus {\mathfrak g}_{\alpha_1}
\cong {\mathfrak s}{\mathfrak l}_2({\mathbb R}) \oplus {\mathbb R}\\
{\mathfrak k}_2 & = & {\mathfrak k}_{\alpha_1} \cong {\mathfrak s}{\mathfrak o}_2.
\end{eqnarray*}
Since ${\mathfrak k}_2$ is one-dimensional, ${\mathfrak v}$ must be a $2$-dimensional linear subspace
of ${\mathfrak n}_1^1$. Let ${\mathfrak v}$ be a ${\mathfrak k}_2$-invariant subspace of ${\mathfrak n}_2^1$.
In order to get a cohomogeneity one action, the normalizer $N_{L_2}({\mathfrak n}_2^1 \ominus {\mathfrak v})$
must act transitively on $F_2 = L_2 \cdot o \cong {\mathbb R}H^2 \times {\mathbb E}$. The only subgroups
of $SL_2({\mathbb R})$ acting transitively on ${\mathbb R}H^2$ are $SL_2({\mathbb R})$ itself and
the parabolic subgroups of $SL_2({\mathbb R})$. However,
${\mathfrak m}_2 \cong {\mathfrak s}{\mathfrak l}_2({\mathbb R})$ acts
irreducibly on ${\mathfrak n}_2^1$, and hence $N_{L_2}({\mathfrak n}_2^1 \ominus {\mathfrak v})$ cannot
be equal to $SL_2({\mathbb R})$. Since $K_2$ is compact and normalizes ${\mathfrak v}$, it also
normalizes ${\mathfrak n}_2^1 \ominus {\mathfrak v}$. If a parabolic subgroup of $SL_2({\mathbb R})$
would normalize ${\mathfrak n}_2^1 \ominus {\mathfrak v}$, then the entire group $SL_2({\mathbb R})$
would normalize ${\mathfrak n}_2^1 \ominus {\mathfrak v}$, which cannot happen. We thus conclude that
$N_{L_2}({\mathfrak n}_2^1 \ominus {\mathfrak v})$ cannot act transitively on $F_2$. This implies that there
is no cohomogeneity one action on $M$ which can be constructed from the choice of $\Phi_2$.
\end{proof}


\end{document}